\documentclass[11pt]{amsart}
\usepackage{geometry}                
\usepackage{graphicx}
\usepackage{epstopdf}
\usepackage{stmaryrd}

\usepackage{hyperref}

\usepackage{hyperref}
\usepackage{stmaryrd}
\usepackage{amsmath}
\usepackage{amssymb}
\usepackage{graphicx}
\usepackage{stackrel}
\usepackage{color}
\usepackage{enumerate}
\usepackage{epsfig}
\usepackage{algorithm}
\usepackage{amsfonts}
\usepackage{graphicx}
\usepackage{amsopn}

\def\db{\textcolor[rgb]{0,0,0.65}}

\def\bit{\begin{itemize}}   \def\eit{\end{itemize}}
\def\ben{\begin{enumerate}} 
\def\een{\end{enumerate}}
\def\benA{\begin{enumerate}[\rm(A)]} 
\def\bena{\begin{enumerate}[\rm(a)]} 
\def\beni{\begin{enumerate}[\rm(i)]}

\def\beq{\begin{equation}}  \def\eeq{\end{equation}}
\def\beqn{\begin{eqnarray}} \def\eeqn{\end{eqnarray}}
\def\beqnn{\begin{eqnarray*}}   \def\eeqnn{\end{eqnarray*}}
\def\barr{\begin{array}}    \def\earr{\end{array}}
\def\bca{\begin{cases}}    \def\eca{\end{cases}}
\def\bfig{\begin{figure}}   \def\efig{\end{figure}}
\def\bpic{\begin{picture}}  \def\epic{\end{picture}}
\def\btab{\begin{tabular}}  \def\etab{\end{tabular}}

\def\sg{{subgradient}}
\def\sd{{subdifferential}}
\def\sda{{subdifferentiable}}
\def\ff{\theta}
\def\gg{\varrho}
\def\hh{h}

\def\bc{\begin{center}}     \def\ec{\end{center}}

\def\qu{\quad}

\def\ms{\medskip}

\def\RR{\mathbb R}
\def\NN{\mathbb N}

\def\RRinf{\RR\cup\{+\infty\} }

\newtheorem{theorem}{\sc Theorem}
\newtheorem{proposition}{\sc Proposition}
\newtheorem{definition}{\sc Definition}
\newtheorem{lemma}{\sc Lemma}
\newtheorem{corollary}{\sc Corollary}
\newtheorem{remark}{\sc Remark}
\newtheorem{example}{\sc Example}
\newtheorem{fact}{\textsc{Fact}}

\def\BT{\begin{theorem}}      \def\ET{\end{theorem}}
\def\BP{\begin{proposition}}  \def\EP{\end{proposition}}
\def\BD{\begin{definition}}   \def\ED{\end{definition}}
\def\BL{\begin{lemma}}        \def\EL{\end{lemma}}
\def\BC{\begin{corollary}}    \def\EC{\end{corollary}}
\def\BR{\begin{remark}\rm }       \def\ER{\end{remark}} 
\def\BE{\begin{example}\rm }      \def\EE{\end{example}}
\def\BF{\begin{fact}}          \def\EF{\end{fact}}
\def\endproof{\hfill $\mb{\small$\Box$}$\ms}
\def\proof{\par\noindent{\textit{Proof}}. \ignorespaces}

\def\cor{Corollary~\ref}
\def\defi{Definition~\ref}

\def\thm{Theorem~\ref}
\def\prop{Proposition~\ref}
\def\lem{Lemma~\ref}

\def\prox{\mathrm{prox}}
\newcommand\dom[1]{\operatorname{dom}(#1)}
\newcommand\Ima[1]{\operatorname{Im}(#1)}

\newcommand\intr[1]{\operatorname{int}(#1)} 
\newcommand\set[1]{{\llbracket 1,#1 \rrbracket}}

\def\lab{\label}

\def\ed{\end{document}}
\def\eq{\eqref}
\def\lp{\left}  \def\rp{\right}
\def\<{\langle} \def\>{\rangle}
\def\leq{\leqslant} \def\geq{\geqslant}
\def\void{\varnothing}
\def\x{\times}
\def\all{\forall\;}

\def\#{\,\sharp\,}

\def\supp{\mathrm{supp}}
\def\d{\partial}

\def\diag{\mathrm{diag}}

\def\wt{\widetilde}

\def\bd{\d}

\def\mb{\mbox}

 \def\B{\mathcal{B}} \def\C{{\mathcal C}} 
   \def\H{{\mathcal H}}

   \def\S{\mathcal{S}}

 \def\V{\mathcal{V}} \def\X{\mathcal{X}}
\def\Y{\mathcal{Y}}

\def\la{\lambda} 
  
\def\ph{\varphi}

\title{A characterization of proximity operators}
\author{R{\'e}mi Gribonval and Mila Nikolova}
\thanks{This work and the companion paper \cite{RGMN2018b} are dedicated to the memory of Mila Nikolova, who passed away prematurely in June 2018. Mila dedicated much of her energy to bring the technical content to completion during the spring of 2018. The first author did his best to finalize the papers as Mila would have wished. He should be held responsible for any possible imperfection in the final manuscript.\protect\\
R. Gribonval (remi.gribonval@inria.fr) was with Univ Rennes, Inria, CNRS, IRISA when this work was conducted. He is now with Univ Lyon, Inria, CNRS, ENS de Lyon, UCB Lyon 1, LIP UMR 5668, F-69342, Lyon, France. M. Nikolova, CMLA, CNRS and ENS de Cachan, Universit{\'e} Paris-Saclay, 94235 Cachan, France.}

\begin{document}

\maketitle

\begin{abstract}
We characterize proximity operators, that is to say functions that map a vector to a solution of a penalized least squares optimization problem. Proximity operators of convex penalties have been widely studied and fully characterized by Moreau.
They are also widely used in practice with nonconvex penalties such as the $\ell^0$ pseudo-norm, yet the extension of Moreau's characterization to this setting seemed to be a missing element of the literature. 
We characterize proximity operators of (convex or nonconvex) penalties as functions that are the subdifferential of some convex potential.
This is proved as a consequence of a more general characterization of so-called Bregman proximity operators
of possibly nonconvex penalties in terms of certain convex potentials. 
As a side effect of our analysis, we obtain a test to verify whether a given function is the proximity operator of some penalty, or not.
Many well-known shrinkage operators are indeed confirmed to be proximity operators. However, we prove that  windowed Group-LASSO and persistent empirical Wiener shrinkage -- two forms of so-called social sparsity shrinkage--
are generally {\em not} the proximity operator of any penalty; the exception is when they are simply weighted versions of group-sparse shrinkage with non-overlapping groups.\protect
\end{abstract}

\section{Introduction and overview} \lab{sec:overview}

Proximity operators have become an important ingredient of nonsmooth optimization, where a huge body of work has demonstrated the power of iterative proximal algorithms to address large-scale variational optimization problems. 
While these techniques have been thoroughly analyzed and understood for proximity operators involving convex penalties, there is a definite trend towards the use of proximity operators of nonconvex penalties such as the $\ell^{0}$ penalty \cite{Blumensath2009,Bredies:2014bh}. 

This paper extends existing characterizations of proximity operators -- which are specialized for convex penalties -- to the nonconvex case. A particular motivation is to understand whether certain thresholding rules known as {\em social sparsity shrinkage}, which have been successfully exploited in the context of certain linear inverse problems, are proximity operators. Another motivation is to characterize when Bayesian estimation with the conditional mean estimator (also known as minimum mean square error estimation or MMSE) can be expressed as a proximity operator. This is the object of a companion paper \cite{RGMN2018b} characterizing when certain variational approaches to address inverse problems can in fact be considered as Bayesian approaches. 

\subsection{Characterization of proximity operators}
Let $\H$ be a Hilbert space equipped with an inner product $\<\cdot,\cdot\>$ and a norm $\|\cdot\|$. 
This includes the case $\H = \RR^n$, and most of the text can be read with this simpler setting in mind.
The proximity operator of a function $\ph: 
\H
\to \RR$ 
maps each $y \in \H$ to the solutions of a penalized least-squares problem
\[ 
y \mapsto \prox_{\ph}(y) := \arg\min_{x \in \H}\lp\{ \tfrac12 \|y-x\|^{2}+\ph(x) \rp\}
\] 
Formally, a proximity operator is set-valued as there may be several solutions to this problem, or the set of solutions may be empty. 
A primary example is the soft-thresholding function $f(y) := y \max(1-1/|y|,0)$,
$y \in \RR$ , which is the proximity operator of the absolute value function $\ph(x):=|x|$. 

Proximity operators can be defined for certain generalized functions $\ph: \H \to \RR \cup \{+\infty\}$. 
A particular example is the projection onto a given closed convex set $\C \subset \H$, which can be written as $\mathrm{proj}_{\C} = \prox_\ph$ with $\ph$ the indicator function of $\C$, i.e., $\ph(x) = 0$ if $x \in \C$, $\ph(x) = +\infty$ otherwise. 
For the sake of precision and brevity, we use the following definition:
\BD\lab{defprox}
Let $\Y \subset \H$ be non-empty. A function $f: \Y \to \H$ is {\em a} proximity operator of  a function $\ph: \H \to \RR \cup \{+\infty\}$ if, and only if,  $f(y) \in \prox_\ph(y)$ for each $y \in \Y$.
\ED
In convex analysis, this corresponds to the notion of a {\em selection} of the set-valued mapping $\prox_{\ph}$.

A characterization of proximity operators {\em of convex lower semicontinuous (l.s.c.) functions} is due to Moreau. It involves the subdifferential $\d\ff(x)$ of a convex l.s.c. function $\ff$ at $x$, i.e., the set of all its subgradients at $x$ \cite[Chapter III.2]{Ekeland83}\footnote{See Section~\ref{sec:notations} for detailed notations and reminders on convex analysis and differentiability in Hilbert spaces.}.\\

\BP~{\cite[Corollary 10.c]{Moreau65}}\label{MoreauNE}
A function $f: \H \to \H$ defined everywhere is the proximity operator of a {\em proper convex l.s.c.} function $\ph : \H \to \RR \cup \{+\infty\}$ if, and only if the following conditions hold jointly:
\bena
\item  there exists a (convex l.s.c.) function $\psi$ such that for each $y \in \H$, $f(y) \in \partial \psi(y)$;
\item $f$ is nonexpansive, i.e.
\[
\|f(y)-f(y')\| \leq \|y-y'\|,\quad \forall y,y' \in \H.
\]
\een
\EP
We extend Moreau's result to possibly nonconvex functions $\ph$ on subdomains of $\H$ by simply relaxing the non-expansivity condition:
\BT\lab{ProCon}
Let $\Y \subset \H$ be non-empty. A function $f: \Y\to \H$ is a proximity operator of a function $\ph: \H \to \RR \cup \{+\infty\}$ 
if, and only if, there exists a convex l.s.c. function $\psi: \H \to \RR \cup \{+\infty\}$ such that for each $y \in \Y$, $f(y) \in \d \psi(y)$. 
\ET
This is proved in Section~\ref{sec:prox} as a particular consequence of our main result, \thm{GProxCh}, which characterizes
 functions such that $f(y) \in \arg\min_{x \in \H} \{ D(x,y) + \ph(x)\}$ for certain types of  data-fidelity terms $D(x,y)$. 
 Among others, the data-fidelity terms covered by \thm{GProxCh} include:
\bit
\item the Euclidean distance $D(x,y) = \tfrac{1}{2}\|y-x\|^2$, which is the data-fidelity associated to proximity operators;
\item its variant $D(x,y) = \tfrac{1}{2} \|y-Mx\|^2$ with $M$ some linear operator; and
\item Bregman divergences \cite{Bregman67}, leading to an analog of \thm{ProCon} to characterize so-called Bregman proximity operators \cite{Censor92} (see \cor{Breprox} in Section~\ref{sec:prox}).
\eit

\thm{GProxCh} further implies that  
the functions 
$\ph$ and $\psi$ in \thm{ProCon} can be chosen such that
\beq
\lab{psivsphi}
\psi(y) = \<y,f(y)\> -\tfrac{1}{2}\|f(y)\|^2-\ph(f(y))
,\qu \all y \in \Y.
\eeq
This is a particular instance of a more general result valid for all considered data-fidelity terms.
Another consequence of \thm{GProxCh} (see \cor{co} in Section~\ref{sec:prox})
is that for the considered data-fidelity terms $D(x,y)$, if $f: \Y \to \H$ can be written as $f(y) \in \arg\min_{x \in \H} \{ D(x,y) + \ph(x)\}$ for some (possibly nonconvex) function $\ph$ and if its image $\Ima{f} := f(\Y)$ is a convex set (e.g., if $\Ima{f}=\H$) then \begin{center}
 {\em the function $x \mapsto D(x,y)+\ph(x)$ is convex on $\Ima{f}$}.
 \end{center}
 This is reminiscent of observations on convex optimization with nonconvex penalties \cite{Parekh:2015wq,Selesnick:2017wa} and on the hidden convexity of conditional mean estimation under additive Gaussian noise \cite{GRIBONVAL:2010:INRIA-00486840:1,NIPS2013_4868,Louchet:2013hs,Advani:2016wg}.
The latter is extended to other noise models in the companion paper \cite{RGMN2018b}.

\subsection{The case of smooth proximity operators}

The smoothness of a proximity operator $f = \prox_\ph$ and that of the corresponding functions $\ph$ and $\psi$, cf \eq{psivsphi}, are inter-related, leading to a characterization of \emph{continuous} proximity operators\footnote{See Section~\ref{sec:notations} for brief reminders on the notion of continuity / differentiability in Hilbert spaces.}.

\BC\label{PCDS}
Let $\Y \subset \H$ be non-empty and open and $f: \Y \to \H$ be $C^{0}$. 
The following are equivalent:
\bena 
\item $f$ is a proximity operator of a function $\ph:\H \to \RR \cup \{+\infty\}$; 
\item  there exists a convex $C^1(\Y)$ 
function $\psi$ such that $f(y) = \nabla \psi(y)$  for each $y \in \Y$.
\een 
\EC
This is established in Section~\ref{localsmooth} as a particular consequence of our second main result, \cor{SmoProx}. There, we also prove that when $f$ is a proximity operator of some $\ph$, the Lipschitz property of $f$ with Lipschitz constant $L$ is equivalent to the convexity of $x \mapsto \ph(x)+\left(1-\tfrac{1}{L}\right)\tfrac{\|x\|^{2}}{2}$. Moreau's characterization (\prop{MoreauNE}) corresponds to the special case $L=1$.
Next, we characterize $C^1$ proximity operators on convex domains more explicitly using the differential of $f$.
\BT\label{jac1}
Let $\Y \subset \H$ be non-empty, open and convex, and $f : \Y \to \H$ be $C^1$.
The following properties are equivalent:
\bena
\item\label{jac0} $f$ is a proximity operator of a function $\ph:\H \to \RR \cup \{+\infty\}$;
\item\label{jaca} there exists a convex $C^{2}(\Y)$
function $\psi$ 
such that $f(y) = \nabla \psi(y)$ for each $y \in \Y$;
\item \label{jacb} the differential $Df(y)$ is a {\em symmetric positive semi-definite} operator\footnote{A continuous linear operator $L: \H \to \H$ is symmetric if $\<x, Ly\> = \<Lx, y\>$ for each $x,y \in \H$. A symmetric continuous linear operator is positive semi-definite if $\<x,Lx\> \geq 0$ for each $x \in \H$. This is denoted $L \succeq 0$. It is positive definite if $\<x,Lx\> >0$ for each nonzero $x \in \H$. This is denoted $L \succ 0$.} for each $y \in \Y$.
\een
\ET
\proof
Since $f$ is $C^{1}$, the equivalence \eq{jac0} $\Leftrightarrow$ \eq{jaca} is a consequence of \cor{PCDS}. 
We now establish \eq{jaca} $\Leftrightarrow$ \eq{jacb}. Since $\Y$ is convex it is simply connected, and as $\Y$ is open by Poincaré's lemma (see \cite[Theorem 6.6.3]{Galbis:2012da} when $\H = \RR^{n}$)  the differential $Df$ is symmetric if, and only if, $f$ is the gradient of some $C^{2}$ function $\psi$. By \cite[Proposition 17.7]{Bauschke:2017ki}, the function $\psi$ is convex iff $\nabla^{2}\psi \succeq 0$ on $\Y$, i.e. iff $Df \succeq 0$ on $\Y$.
\endproof

\BC \lab{jac2}
Let $\Y \subset \H$ be a set with non-empty interior $\intr{\Y} \neq \void$, $y \in \intr{\Y}$,
 and $f: \Y \to \H$ be a proximity operator. If $f$ is $C^1$ in a neighborhood of 
$y$,
 then $Df(y)$ is symmetric positive semi-definite. 
\EC
\proof Restrict $f$ to any open convex neighborhood $\Y' \subset \Y$ of $y$ and apply \thm{jac1}.
\endproof

\BR
Differentials are perhaps more familiar to some readers in the context of multivariate calculus: when $y = (y_{j})_{j=1}^{n} \in \H = \RR^{n}$ and $f(y) = (f_{i}(y))_{i=1}^{n}$, $Df(y)$ is identified to the Jacobian matrix 
\[
Jf(y) = (\tfrac{\partial f_{i}}{\partial y_{j}})_{1 \leq i,j \leq n}.
\]
The rows of $Jf(y)$ are the transposed gradients $\nabla f_{i}(y)$. The differential is symmetric if the mixed derivatives satisfy $\tfrac{\partial f_{i}}{\partial y_{j}} = \tfrac{\partial f_{j}}{\partial y_{i}}$ for all $i \neq j$. When $n=3$, this corresponds to $f$ being an {\em irrotational vector field}. More generally, this characterizes the fact that $f$ is a so-called {\em conservative field}, i.e., a vector field that is the gradient of some potential function. As the Jacobian is the Hessian of this potential, it is positive definite if the potential is convex.
\ER

Finally we provide conditions ensuring that $f$ is a proximity operator and that $f(y)$ is the only critical point of the corresponding optimization problem. 
\BC\lab{uniqueC1}
Let $\Y \subset \H$ be open and convex, and $f: \Y \to \H$ be $C^1$ with $Df(y) \succ 0$ on $\Y$. 
Then $f$ is injective and there is $\ph: \H \to \RR \cup \{+\infty\}$ such that $\prox_\ph(y) = \{f(y)\},\ \all y \in \Y$ and $\dom \ph = \Ima{f}$. 
Moreover, if $Df(y)$ is {\em boundedly invertible} on $\Y$ then $\ph$ is $C^{1}$ on $\Y$ and for each $y \in \Y$, the only critical point of $x \mapsto \tfrac{1}{2}\|y-x\|^2+\ph(x)$
is $x=f(y)$.
\EC
This is established in Appendix~\ref{pfuniqueC1}.
\BR \lab{rem:bndinv}
In finite dimension $\H = \RR^{n}$, $Df(y)$ is boundedly invertible as soon as $Df(y) \succ 0$, hence we only need to assume that $Df(y) \succ 0$ to conclude that $f(y)$ is the unique critical point. This is no longer the case in infinite dimension. Indeed, consider $\H = \ell^{2}(\NN)$ and $f: y = (y_{n})_{n \in \NN} \mapsto f(y) := (y_{n}/(n+1))_{n \in \NN}$. As $f$ is linear, its differential is $Df(y) = f$ for every $y \in \H$. As $\langle f(y),y\rangle = \sum_{n \in \NN} y_{n^{2}}/(n+1) > 0$ for each nonzero $y \in \H$ we have $Df(y) \succ 0$ but its inverse is unbounded. Given $n \in \NN$ and $z \in \RR$ we have $z/(n+1) = \arg\min_{x \in \RR} \tfrac{1}{2}(z-x)^{2}+nx^{2}/2$, hence $f = \prox_{\ph_{0}}$ with $\ph_{0}: x = (x_{n})_{n \in \NN}  \mapsto \ph_{0}(x) := \sum_{n \in \NN} n x_{n}^{2}/2$. Setting $\ph(x) = \ph_{0}(x)$ for $x \in \Ima{f}$, $\ph(x) = +\infty$ otherwise, we have $\prox_{\ph} = f$ and $\dom{\ph} = \Ima{f}  = \{x \in \H, \sum_{n \in \NN} (n+1)^{2}x_{n}^{2} < \infty\}$. Yet, as no point in $\dom{\ph}$ admits any open neighborhood in $\H$, $\ph$ is nowhere differentiable and {\em every} $x \in \H$ is a critical point of $x \mapsto \tfrac{1}{2}\|y-x\|^{2}+\ph(x)$.
\ER

{\bf Terminology.} Proximity operators often appear in the context of penalized least squares regression, where $\ph$ is called a {\em penalty}, and from now on we will adopt this terminology. In light of \cor{PCDS}, a continuous proximity operator is exactly characterized as a gradient of a convex function $\psi$. In the terminology of physics, a proximity operator is thus a {\em conservative field} associated to a {\em convex potential}. In the language of convex analysis, subdifferentials of convex functions are characterized as maximal cyclically monotone operators \cite[Theorem B]{Rockafellar:1970fu}.

\subsection{Organization of the paper}

The proof of our most general results, \thm{GProxCh} and \cor{SmoProx} (and the fact that they imply \thm{ProCon}, \eqref{psivsphi}, \cor{co} and \cor{PCDS}) are established in Section~\ref{sec:prox}, where we also discuss their consequences in terms of Bregman proximity operators and illustrate them on concrete examples. As \thm{ProCon} and its corollaries characterize whether a function $f$ is a proximity operator and study its smoothness in relation to that of the corresponding penalty and potential, they are particularly useful when $f$ is not  {\em explicitly built} as a proximity operator. This is the case of so-called social shrinkage operators (see e.g. \cite{Kowalski:2013vx}). We conclude the paper by showing  in Section~\ref{sec:social} that social shrinkage operators are generally {\em not} the proximity operator of any penalty.

\subsection{Discussion}

In light of the extension to nonconvex penalties of Moreau's characterization of proximity operators of convex (l.s.c.) penalties (\prop{MoreauNE}), the nonexpansivity of the proximity operator $f$ determines whether the  underlying penalty $\ph$ is convex or not. While non-expansivity certainly plays a role in the convergence analysis of iterative proximal algorithms based on convex penalties, the adaptation of such an analysis when the proximity operator is Lipschitz rather than nonexpansive, using \prop{MoreauNE}, is an interesting perspective. 

The characterization of smooth proximity operators as the gradients of convex potentials, which also appear in optimal transport (see e.g., \cite{Villani:2008vs}), suggests that further work is needed to better understand the connections between these concepts and tools. This could possibly lead to simplified arguments where the strong machinery of convex analysis may be used more explicitly despite the apparent lack of convexity of the optimization problems associated to nonconvex penalties.

\section{Main results}\lab{sec:prox}
We now state our main results, \thm{GProxCh} and \cor{SmoProx}, and prove a number of their consequences including \thm{ProCon}, \eqref{psivsphi}, \cor{co} and \cor{PCDS} which were advertized in Section~\ref{sec:overview}.
The most technical proofs are postponed to the Appendix.

\subsection{Detailed notations} \label{sec:notations}
The indicator function of a set $\S$ is denoted
\[ \chi_{\S}(x):=\lp\{ \barr{lll} 0&\mb{if}&x\in \S, \\ +\infty&\mb{if}& x\not\in \S.\earr\right.\]
The domain of a function $\ff : \H \to \RR \cup \{+\infty\}$ is defined and denoted by $\dom{\ff} := \{x \in \H\mid \ff(x) < \infty\}$. 
Given $\Y \subset \H$ and a function $f: \Y \to \H$, the image of $\Y$ under $f$ is denoted by $\Ima{ f }$.
A function $\ff : \H \to \RR \cup \{+\infty\}$ is proper iff there is $x\in\H$ such that $\ff(x) < +\infty$, i.e., $\dom{\ff} \neq \void$. It is lower semicontinuous (l.s.c.) if for each $x_0 \in \H$, $\liminf_{x \to x_0} \ff(x) \geq \ff(x_0)$, or equivalently if the set $\{x \in \H: \ff(x)>\alpha\}$ is open for every $\alpha \in \RR$. 
A subgradient of a convex function $\ff: \H \to \RR \cup \{+\infty\}$ at $x$ is any $u \in \H$ such that $
\ff(x')-\ff(x) \geq \< u,x'-x\>, \forall x' \in \H$. A function with $k$ continuous derivatives\footnote{see  Appendix~\ref{frechet} for some reminders on Fr{\'e}chet derivatives in Hilbert spaces.} is called a $C^k$ function. The notation $C^k(\X)$ is used to specify a $C^k$ function on an open domain $\X$. Thus $C^0$ is the space of continuous functions, whereas $C^1$ is the space of continuously differentiable functions  \cite[p. 327]{Cartan77}. The gradient of a $C^{1}$ scalar function $\ff$ at $x$ is denoted $\nabla \ff(x)$.

The segment between two elements $x,x' \in \H$ is the set $[x,x'] := \{tx+(1-t)x', t \in [0,1]\}$. A finite union of segments $[x_{i-1},x_{i}]$, $1 \leq i \leq n$, $n \in \NN$, where $x_{0}=x$ and $x_{n}=x'$ is called a polygonal path between $x$ and $x'$. A non-empty subset $\C \subset \H$ is polygonally connected iff between each pair $x,x' \in \C$ there is a polygonal path with all its segments included in $\C$, $[x_{i-1},x_{i}] \subset \C$.
\BR
The notion of polygonal-connectedness is a bit stronger than that of connectedness. Indeed, polygonal-connectedness implies the classical topological property of path-connectedness, which in turn implies connectedness. However there are path-connected sets that are not polygonally-connected -- e.g., the unit circle in $\RR^2$ is path-connected, but no two points are polygonally-connected, and there are connected sets that are not path-connected. Yet, every {\em open} connected set is polygonally-connected, see \cite[Theorem 2.5.2]{Galbis:2012da} for a statement in $\RR^{n}$.
\ER

\subsection{Main theorem}

\BT \lab{GProxCh}
Consider $\H$ and $\H'$ two Hilbert spaces\footnote{For the sake of simplicity we use the same notation $\<\cdot,\cdot\>$ for the inner products $\<x,A(y)\>$ (between elements of $\H$) and $\<B(x),y\>$ (between elements of $\H'$). The reader can inspect the proof of \thm{GProxCh} to check that the result still holds if we consider {\em Banach spaces} $\H$ and $\H'$, $\H^\star$ and $(\H')^\star$ their duals, and $A: \Y \to \H^\star$, $B: \H \to (\H')^\star$.}, and $\Y \subset \H'$ a non-empty set.
 Let $a: \Y \to \RR \cup \{+\infty\}$, $b: \H \to \RR \cup \{+\infty\}$, $A: \Y \to \H$ and $B: \H \to \H'$ be arbitrary functions. 
 Consider $f: \Y \to \H$  and denote $\Ima{f}$ the image of $\Y$ under $f$.
\bena
\item \lab{pro1} Let $D(x,y) := a(y)-\<x,A(y)\>+b(x)$. The following properties are equivalent:
\beni
\item \lab{proa} there is $\ph: \H \to \RR \cup \{+\infty\}$ such that
$f(y) \in \arg\min_{x \in \H} \{D(x,y)+\ph(x)\}$ for each $y \in \Y$;
\item \lab{prob} there is a convex l.s.c. $g: \H \to \RR \cup \{+\infty\}$ such that $A(f^{-1}(x))\subset \d g \lp( x \rp)$ for each $x \in \Ima{f}$;
 \een
When they hold, $\ph$ (resp. $g$) can be chosen given $g$ (resp. $\ph$) so that $g(x)+\chi_{\Ima{f}} = b(x)+\ph(x)$.
 \item\lab{pro2} Let $\ph$ and $g$ satisfy \eq{proa} and \eq{prob}, respectively, and let $\C \subset \Ima{f}$ be
polygonally connected. Then there is $K \in \RR$ such that
 \beqn
 \lab{ggg} g(x) &=& b(x)+\ph(x) + K,\qu \all x \in \C.
 \eeqn 
 \item\lab{pro3} Let $\wt D(x,y) := a(y)-\<B(x),y\>+b(x)$.  The following properties are equivalent:
\beni
\item \lab{proabis} there is $\ph: \H \to \RR \cup \{+\infty\}$ such that 
$f(y) \in \arg\min_{x \in \H} \{\wt D(x,y)+\ph(x)\}$ for each $y \in \Y$;
\item \lab{proc} there is a convex l.s.c. $\psi: \H' \to \RR \cup \{+\infty\}$ such that $B(f(y))\in \d \psi(y)$ for each $y \in \Y$.
 \een 
$\ph$ (resp. $\psi$) can be chosen given $\psi$ (resp. $\ph$) so that $\psi(y) = \langle B(f(y),y\rangle-b(f(y))-\ph(f(y))$ on $\Y$.
 \item\lab{pro4} Let $\ph$ and $\psi$ satisfy \eq{proabis} and \eq{proc}, respectively, and let $\C' \subset \Y$ be
 polygonally connected. Then there is $K' \in \RR$ such that
 \beqn
 \lab{psi} \psi(y) &=& \<B(f(y)),y\>-b(f(y))-\ph(f(y)) + K',\qu \all y \in \C'.
 \eeqn 
 \een
\ET

The proof of \thm{GProxCh} is postponed to Appendix~\ref{pf:GProxCh}. As stated in \eq{pro1} (resp. \eq{pro3}), the functions can be chosen such that the relation \eq{ggg} (resp. \eq{psi}) holds on $\Ima{f}$ (resp. on $\Y$) with $K=K'=0$. As the functions $\ph,g,\psi$ are at best defined up to an additive constant, we provide in \eq{pro2} (resp. \eq{pro4}) conditions ensuring that adding a constant is indeed the unique degree of freedom. The role of polygonal-connectedness will be illustrated on examples in Section~\ref{sec:examplessocshrink}.
\BE \label{lininv}
In the context of linear inverse problems one often encounters optimization problems involving functions expressed as $\tfrac{1}{2}\|y-Mx\|^2 + \ph(x)$ with $M$ some linear operator. Such functions fit into the framework of \thm{GProxCh} using $a(y):=\tfrac{1}{2}\|y\|^2$, $b(x):=\tfrac{1}{2}\|Mx\|^2$, $A(y):=M^\star y$, and $B(x) := Mx$,  where $M^\star$ is the adjoint of $M$.
Among other consequences one gets that $f:\Y \to \H$ is a generalized proximity operator of this type for some penalty $\varphi$ if, and only if, there is a convex l.s.c. $\psi$ such that $M f(y) \in \d \psi(y)$ for each $y \in \Y$.
\EE
Examples where the data-fidelity term is a so-called Bregman divergence are detailed in Section~\ref{bred} below. This covers the case of standard proximity operators where  $D(x,y) = \tfrac{1}{2}\|y-x\|^{2}$.

\subsection{Convexity in proximity operators of nonconvex penalties}
\lab{cvxnncvx}
An interesting consequence of \thm{GProxCh} is that the optimization problem associated to (generalized) proximity operators is in a sense always convex, even when the considered penalty $\ph$ is not convex.

\BC\lab{co}
Consider $\H,\H'$ two Hilbert spaces.
Let $\Y \subset \H'$ be non-empty and $f:\Y\to\H$.
Assume that  there is $\ph: \H \to \RR \cup \{+\infty\}$ such that $f(y) \in \arg\min_{x \in \H} \lp\{D(x,y)+\ph(x)\rp\}$ for each $y \in \Y$,  with $D(x,y) = a(y)-\<x,A(y)\>+b(x)$ as in \thm{GProxCh}\eq{pro1}.
Then
\bena
\item \lab{coa} the function $x \mapsto b(x)+\ph(x)$ is convex on each convex subset $\C\subset\Ima{f}$;
\item \lab{cob}  if $\Ima{f}$ is convex, then the function $x \in \Ima{f} \mapsto D(x,y)+\ph(x)$ is convex, $\all y\in \Y$.
\een
Similarly, if there is $\ph: \H \to \RR \cup \{+\infty\}$ such that $f(y) \in \arg\min_{x \in \H} \lp\{\wt D(x,y)+\ph(x)\rp\}$ for each $y \in \Y$,  with $\wt D(x,y) = a(y)-\<B(x),y\>+b(x)$ as in \thm{GProxCh}\eq{pro3} then
 $y \mapsto \<B(f(y)),y\>-b(f(y))-\ph(f(y))$ is convex on each convex subset $\C'\subset\Y$.
\EC
\proof 
\eq{coa} follows from \thm{GProxCh}\eq{pro1}-\eq{pro2}.  \eq{cob}
follows from \eq{coa} and the definition of $D$. The proof of the result with $\wt D$ instead of $D$ is similar.
\endproof

\cor{co}\eq{cob} might seem surprising as, given a nonconvex penalty $\ph$, one may expect the optimization problem $\min_x D(x,y)+\ph(x)$ to be nonconvex.
However, as noticed e.g. by \cite{Nikolova:kk,Parekh:2015wq,Selesnick:2017wa}, there are nonconvex penalties such that this problem with $D(x,y) := \tfrac{1}{2}\|y-x\|^2$ is in fact convex. \cor{co} establishes that this convexity property indeed holds whenever the image $\Ima{f}$ of the resulting function $f$ is a convex set.
A particular case is that of functions $f$ built as conditional expectations in the context of additive Gaussian denoising, which have been shown \cite{GRIBONVAL:2010:INRIA-00486840:1} to be proximity operators. Extensions of this phenomenon for conditional mean estimation with other noise models are discussed in the companion paper \cite{RGMN2018b}.

\subsection{Application to Bregman proximity operators} \label{bred}

The squared Euclidean norm is a particular {\em Bregman divergence}, and \thm{GProxCh} characterizes generalized proximity operators defined with such divergences. The Bregman divergence, known also as $D$-function, 
was introduced in \cite{Bregman67} for strictly convex differentiable functions on so-called linear topological spaces. 
For the goals of our study, it will be enough to 
consider that $\hh:\H\to\RRinf$ is proper, convex and differentiable on a Hilbert space.

\BD \lab{bregman}
Let $\hh:\H\to\RRinf$ be proper convex and differentiable on its open domain $\dom{\hh}$. 
The Bregman divergence (associated with $\hh$) between $x$ and $y$  
is defined by
\beq\label{->} 
 D_\hh :\H\x\H \to [0,+\infty]:(x,y)\to
\begin{cases}
\hh(x)-\hh(y) -\<\nabla\hh(y),x-y\>, & \text{if}\ y \in \dom{\hh};\\
+\infty, & \text{otherwise}
\end{cases}
\eeq
\ED

In \thm{GProxCh}\eq{pro1} one obtains $D(x,y) = D_\hh(x,y)$ by setting $a(y) = +\infty$ and $A(y)$ arbitrary if $y \notin \dom{\hh}$
and, for $y \in \dom{\hh}$ and each $x \in \H$,
\beq\label{f->} a(y):=\<\nabla \hh(y),y\> - \hh(y)\qu b(x):=\hh(x)\qu \qu \mb{and}
\qu A(y)=\nabla \hh (y)\eeq
The lack of symmetry of the Bregman divergence suggests to  consider also $D_h(y,x)$. 
In \thm{GProxCh}\eq{pro3} one obtains $\wt D(x,y) = D_h(y,x)$ using $b(x) = +\infty$ and $B(x)$ arbitrary for $x \notin \dom{\hh}$ and, for $x \in \dom{\hh}$ and each $y \in \H$,
\beq\label{f<-} a(y):=\hh(y)\qu\qu b(x):=\<\nabla\hh(x),x\>-\hh(x)\qu\mb{and}\qu B(x)=\nabla\hh(x)\eeq

The next claim is an application of \thm{GProxCh}
with $D(x,y) = D_\hh(x,y)$ and $\wt D(x,y) = D_\hh(y,x)$.
We thus consider the so-called Bregman proximity operators which were introduced in \cite{Censor92}. 
We will focus on the characterization of 
these operators defined by $y\mapsto \arg\min_{x \in \H} \{D_\hh(x,y)+\ph(x)\}$ and $y\mapsto \arg\min_{x \in \H} \{D_\hh(y,x)+\ph(x)\}$. Such operators have been further studied in \cite{Bauschke:2003jf} with an emphasis on the notion of viability, which is essential for these operators to be useful in the context of iterative algorithms.

\BC \lab{Breprox}
Consider $f: \Y \to \H$. Let $\hh:\H\to\RRinf$ be a proper convex function that is differentiable on its open domain $\dom{\hh}$.
Let $D_h$ read as in \eq{->}.
\bena
\item \lab{Bpro1} The following properties are equivalent:
\beni
\item \lab{Bproa} there is $\ph: \H \to \RR \cup \{+\infty\}$ such that
$f(y) \in \arg\min_{x \in \H} \{D_\hh (x,y)+\ph(x)\}$, $\all y \in \Y$;
\item \lab{Bprob} there is a convex l.s.c. $g: \H \to \RR \cup \{+\infty\}$ s.t. $\nabla\hh(f^{-1}(x))\subset \d g \lp( x \rp)$,
 $\all x \in \Ima{f}$;
 \een
 When they hold, $\ph$ (resp. $g$) can be chosen given $g$ (resp. $\ph$) so that $g(x)+\chi_{\Ima{f}} = h(x)+\ph(x)$.
\item\lab{Bpro2} Let $\ph$ and $g$ satisfy \eq{Bproa} and \eq{Bprob}, respectively, and let $\C \subset \Ima{f}$ be polygonally connected. Then there is $K \in \RR$ such that
\[ 
  g(x) = \hh(x)+\ph(x) + K,\qu \all x \in \C.
\]  

 \item\lab{Bpro3} The following properties are equivalent:
\beni
\item \lab{Bproabis} there is $\ph: \H \to \RR \cup \{+\infty\}$ such that
$f(y) \in \arg\min_{x \in \H} \{D_\hh(y,x)+\ph(x)\}$, $\all y \in \Y$;
\item \lab{Bproc} there is a convex l.s.c. $\psi: \H \to \RR \cup \{+\infty\}$ such that $\nabla\hh(f(y))\in \d \psi(y)$, $\all y \in \Y$.
 \een 
  $\ph$ can be chosen given $\psi$ (resp. $\psi$ given $\ph$) s.t. $\psi(y) = \langle \nabla h(f(y)),y-f(y)\rangle+h(f(y))-\ph(f(y)), \forall y \in  \Y$.
 \item\lab{Bpro4} Let $\ph$ and $\psi$ satisfy \eq{proabis} \eq{proc}, respectively, and let $\C' \subset \Y$ be polygonally connected. Then there is $K' \in \RR$ such that
\[ 
 \psi(y) =\big\<\nabla\hh(f(y)),y-f(y) \big\>+\hh(f(y))-\ph(f(y)) + K',\qu \all y \in \C'.
\]  
 \een
\EC

\proof \eq{Bpro1} and \eq{Bpro2} use 
\eq{f->}. Further,  \eq{Bpro3} and \eq{Bpro4} use 
\eq{f<-}.
\endproof

\subsection{Specialization to (standard) proximity operators} \label{sec:classicalprox}

Standard (Hilbert space) proximity operators correspond to taking as the Bregman divergence $D_{h}(x,y) = \tfrac{1}{2} \|y-x\|^2$, which is associated to $h(x) := \tfrac{1}{2}\|x\|^{2}$. An immediate consequence of \cor{bred} is the following theorem, which implies \thm{ProCon} and \eq{psivsphi}.

\BT \lab{ProxCh}
Let $\Y\subset \H$ be non-empty, and $f: \Y \to \H$. 
\bena
\item  \lab{tpro0} The following properties are equivalent:
\beni
\item \lab{tproa} there is $\ph: \H \to \RR \cup \{+\infty\}$ such that
$f(y) \in \prox_{\ph}(y)$ for each $y \in \Y$;
\item \lab{tprob} there is a convex l.s.c. $g: \H \to \RR \cup \{+\infty\}$ such that $f^{-1}(x) \subset \d g \lp( x \rp)$ for each $x \in \Ima{f}$;
\item \lab{tproc} there is a convex l.s.c. $\psi: \H \to \RR \cup \{+\infty\}$ such that $f(y) \in \d \psi(y)$ for each $y \in \Y$.
 \een
When they hold, there exists a choice of $\ph,g,\psi$ satisfying \eq{tproa}-\eq{tprob}-\eq{tproc} such that
  \beqnn
  g(x)+\chi_{\Ima{f}} &=& \tfrac{1}{2}\|x\|^{2}+\ph(x),\qu \all x \in \H;\\
  \psi(y) &=& \langle y,f(y)\rangle-\tfrac{1}{2}\|f(y)\|^{2}-\ph(f(y)),\qu \all y \in \Y.
  \eeqnn
 \item \lab{tpro} Let $\ph$, $g$ and $\psi$ satisfy \eq{tproa}, \eq{tprob} and \eq{tproc}, respectively. 
 Let $\C \subset \Ima{f}$ and $\C' \subset \Y$ be polygonally connected. Then there exist $K,K' \in \RR$ such that
 \beqn
 \lab{tggg} g(x) &=& \tfrac{1}{2}\|x\|^2+\ph(x) + K,\qu \all x \in \C;\\
 \lab{tpsi} \psi(y) &=& \<y,f(y)\> 
-\tfrac{1}{2}\|f(y)\|^2-\ph(f(y))
+ K',\qu \all y \in \C'.
 \eeqn
 \een
\ET

\subsection{Local smoothness of proximity operators}
\label{localsmooth}
\thm{ProxCh} characterizes proximity operators in terms of three functions: a (possibly nonconvex) penalty $\ph$, a convex potential $\psi$, and another convex function $g$. As we now show, the 
properties of these functions are tightly inter-related. 
First we extend Moreau's characterization (\prop{MoreauNE}) as follows:
\BP  \lab{MoreauNEext}
Consider $f: \H \to \H$ defined everywhere, and $L>0$. The following are equivalent:
\ben
\item \lab{MNEext1} there is $\ph: \H \to \RR \cup \{+\infty\}$ s.t. $f(y) \in \prox_{\ph}(y)$ on $\H$, and $x \mapsto \ph(x)+(1-\tfrac{1}{L})\tfrac{\|x\|^{2}}{2}$ is convex l.s.c;
\item  the following conditions hold jointly:
\bena
\item \lab{MNEexta} there exists a (convex l.s.c.) function $\psi$ such that for each $y \in \H$, $f(y) \in \partial \psi(y)$;
\item \lab{MNEextb} $f$ is $L$-Lipschitz, i.e.
\[
\|f(y)-f(y')\| \leq L \|y-y'\|,\quad \forall y,y' \in \H.
\]
\een
\een
\EP
\proof \eq{MNEext1}$\Rightarrow$\eq{MNEexta}.
Simply observe that $f$ is a proximity operator and use \thm{ProxCh}\eq{tproa}$\Rightarrow$\eq{tproc}.\\
\eq{MNEext1}$\Rightarrow$\eq{MNEextb}.
The function $\tilde{\ph}(z) := \tfrac{1}{L} (\ph(Lz)+(1-\tfrac{1}{L})\tfrac{\|Lz\|^{2}}{2})$ is convex l.s.c. by assumption. We prove below that $\tilde{f} := f/L$ is a proximity operator of $\tilde{\ph}$. By \prop{MoreauNE} $\tilde{f}$ is thus non-expansive, i.e., $f$ is $L$-Lipschitz.\\
To show $\tilde{f}(y) \in \prox_{\tilde{\ph}}(y)$ for each $y \in \H$, observe that $\ph(x) = L \tilde{\ph}(x/L)-(1-\tfrac{1}{L})\tfrac{\|x\|^{2}}{2}$. For each $x \in \H$
\beqnn
\tfrac{1}{2}\|y-x\|^{2}+\ph(x) 
&=& \tfrac{\|y\|^{2}}{2}-\langle y,x\rangle +  \tfrac{\|x\|^{2}}{2}+L\tilde{\ph}(x/L)-(1-\tfrac{1}{L})\tfrac{\|x\|^{2}}{2}
= \tfrac{\|y\|^{2}}{2}-\langle y,x\rangle +  \tfrac{\|x\|^{2}}{2L}+L\tilde{\ph}(x/L)\\
&=& \tfrac{\|y\|^{2}}{2}-L\langle y,z\rangle + L \tfrac{\|z\|^{2}}{2}+L\tilde{\ph}(z)\\
&=& (1-L)\tfrac{\|y\|^{2}}{2}+L\left(\tfrac{1}{2}\|y-z\|^{2}+\tilde{\ph}(z)\right),\qu \text{with}\ z = x/L.
\eeqnn
Since $x=f(y)$ is a minimizer of the left-hand-side, $z = f(y)/L = \tilde{f}(y)$ is a minimizer of the right hand side, hence $\tilde{f}$ is a proximity operator of $\tilde{\ph}$ as claimed.\\
\eq{MNEexta} and \eq{MNEextb}$\Rightarrow$\eq{MNEext1}.
By \eq{MNEexta} the function $\tilde{\psi}(y) := \psi(y)/L$ is convex l.s.c and $f(y)/L \in \partial \tilde{\psi}(y)$. By \thm{ProxCh}\eq{tproc}$\Rightarrow$\eq{tproa}$\tilde{f} := f/L$ is therefore a proximity operator.  
Since $f$ is $L$-Lipschitz, $\tilde{f}$ is non-expansive hence by \prop{MoreauNE} $\tilde{f}$ is a proximity operator of some {\em convex l.s.c} penalty $\tilde{\ph}$. The function $\ph(x) := L \tilde{\ph}(x/L)-(1-\tfrac{1}{L})\tfrac{\|x\|^{2}}{2}$ is such that  $\ph(x)+(1-\tfrac{1}{L})\tfrac{\|x\|^{2}}{2} = L \tilde{\ph}(x/L)$ is convex l.s.c. as claimed. By the same argument as above, as $z = \tilde{f}(y)$ is a minimizer of $\tfrac{1}{2}\|y-z\|^{2}+\tilde{\ph}(x)$, $x = Lz = f(y)$ is a minimizer of $\tfrac{1}{2}\|y-x\|^{2}+\ph(x)$,  showing that $f$ is indeed a proximity operator of $\ph$. 
\endproof

Next we consider additional properties of these functions.
\BC\label{SmoProx}
Let $\Y \subset \H$ and $f : \Y \to \H$. 
Consider three functions $\ph$, $g$, $\psi$ on $\H$ 
satisfying the equivalent properties \eq{tproa}, \eq{tprob} and 
\eq{tproc} of \thm{ProxCh}, respectively. Let $k \geq 0 $ be an integer. 
\bena
\item \lab{sp1}
Consider an open set $\V \subset \Y$. The following two properties are equivalent:
\beni
\item \lab{spa} $\psi$ is $C^{k+1}(\V)$;
\item \lab{spb} $f$ is 
$C^k(\V)$; 
\een
When one of them holds, we have $f(y) = \nabla \psi(y), \forall y \in \V$.
\item \lab{sp2} Consider an open set $\X \subset \Ima{f}$. The following three properties are equivalent:
\beni
\item \lab{spc} $\ph$ is $C^{k+1}(\X)$;
\item \lab{spd} $g$ is $C^{k+1}(\X)$;
\item \lab{spe} the restriction $\wt f$ of $f$ to the set $f^{-1}(\X)$ is injective and $(\wt f)^{-1}$ is $C^k(\X)$.
\een
When one of them holds, $\wt f$ is a bijection between $f^{-1}(\X)$ and $\X$, and  we have 
\[
(\wt f)^{-1}(x) = \nabla g(x) = x+\nabla\ph (x), \forall x \in \X.
\]
\een
\EC
Before proving this corollary, let us first mention that the characterization of any {\em continuous} proximity operator $f$ as the gradient  of a $C^1$ convex potential $\psi$, i.e., $f = \nabla \psi$, is a direct consequence of \cor{SmoProx}\eqref{sp1} and \thm{ProCon}. This establishes  \cor{PCDS} from Section~\ref{sec:overview}.

The proof of \cor{SmoProx} relies on the following technical lemma which we prove in Appendix~\ref{lsp} as a consequence of \cite[Prop 17.41]{Bauschke:2017ki}.
\BL\label{SdCont}
Consider a function $\gg: \H \to \H$, a function $\ff: \H \to \RR \cup \{+\infty\}$ and an open set $\X\subset \dom{\gg} \cap \dom{\ff} \subset \H$. Assume that $\ff$ is {\sda}  
at each $x \in \X$ and that 
\beq\lab{ginf} \all x \in \X\qu\qu \gg(x) \in \bd \ff(x)\eeq
Then the following statements are equivalent: 
\bena 
\item \lab{SdConta}
$\gg$ is continuous on $\X$;
\item \lab{SdContb} 
$\ff$ is continuously differentiable on $\X$
i.e., its gradient $\nabla \ff(x)$ is continuous on $\X$.  
\een
When one of the statements holds, $\{\gg(x)\} = \{\nabla \ff(x)\}=\bd\ff(x)$ for each $x \in \X$. 
\EL
\proof [Proof of \cor{SmoProx}]

\eq{spa} $\Leftrightarrow$ \eq{spb} By assumption $\psi$ satisfies \thm{ProxCh}\eq{proc}, i.e., $f(y) \in \d \psi(y)$, $\all y \in \V$.
By \lem{SdCont} with $\gg:=f$ and the convex function $\ff:=\psi$, 
$f $ is $C^0(\V)$ if and only if $\psi$ is $C^1(\V)$ and 
when one of these holds, $f = \nabla \psi$ on $\V$. This proves the result for $k=0$. The extension to $k \geq 1$ is trivial. 

\eq{spc} $\Leftrightarrow$ \eq{spd} Consider $x \in \V$. As $\V$ is open there is an open ball $\B_x$ such that $x \in \B_x \subset \V$. Noticing that $\B_x$ is polygonally connected, by \thm{ProxCh}-\eq{tpro}, there is $K \in \RR$ such that $g(x') = \tfrac{1}{2}\|x'\|^2+\ph(x')+K$ for each $x' \in \B_x$. Hence $g$ is $C^{k+1}(\B_x)$ if and only if $\ph$ is $C^{k+1}(\B_x)$, and $\nabla g(x') = x' + \nabla \ph(x')$ on $\B_x$. As this holds for each $x \in \V$, the equivalence holds on $\V$.

\eq{spd} $\Rightarrow$ \eq{spe} 
By \eq{spd}, $g$ is $C^{k+1}(\X)$ hence $\d g(x) = \{\nabla g(x)\}$ for each $x \in \X$. By \thm{ProxCh}\eq{prob}, $f^{-1}(x)  \subset \d g(x)$ for each $x \in \Ima{f}$. Combining both facts yields 
\beq\lab{in} y = \nabla g(f(y))\qu \all y\in f^{-1}(\X).\eeq 
Consider $y,y' \in f^{-1}(\X)$ such that $f(y) = f(y')$. Then $y = \nabla g(f(y)) = \nabla g(f(y')) = y'$, 
which shows that $f$ is injective on $f^{-1}(\X)$. 
Consequently, $\wt f$ is a bijection between $f^{-1}(\X)$ and $\X$, 
hence the inverse function $(\wt f)^{-1} $ is well defined.
Inserting $y=(\wt f)^{-1} (x)$ into \eq{in} yields $(\wt f)^{-1} (x) = \nabla g(x)$ for each $x \in \X$. 
Then, since $g$ is $C^{k+1}(\X)$, it follows that $(\wt f)^{-1} $ is $\C^k(\X)$.

\eq{spe} $\Rightarrow$ \eq{spd} 
Consider $x \in \X$. As  $\wt f$ is injective on $f^{-1}(\X)$ by \eq{spe}, there is a unique $y \in f^{-1}(\X)$ such that $x=f(y)$. Using that $f^{-1}(x) \subset \d g(x)$ by \thm{ProxCh}\eq{prob} shows that $(\wt f)^{-1}(x) = y\in \d g(x)$. Since $(\wt f)^{-1}$ is $C^k(\X)$, using \lem{SdCont}  with $\gg:=(\wt f)^{-1}$ and $\ff:=g$ proves that 
\[ 
(\wt f)^{-1}(x) =\nabla g(x) \qu \all x\in\X
\] 
Since $(\wt f)^{-1}$ is $C^k(\X)$ it follows that $g$ is $C^{k+1}(\X)$. 
\endproof

\subsection{Illustration using classical examples}\lab{sec:examplessocshrink}

\thm{ProCon} and its corollaries characterize whether a function $f$ is a proximity operator. This is particularly useful when $f$ is not  {\em explicitly built} as a proximity operator. We illustrate this with a few examples. We begin with $\H = \RR$, where proximity operators happen to have a particularly simple characterization.

\BC \lab{coscal}
Let $\Y \subset \RR$ be non-empty. A function $f: \Y \to \RR$ is the proximity operator of some penalty $\ph$ if, and only if, $f$ is  nondecreasing. 
\EC
\proof By \thm{ProCon} we just need to prove that a scalar function $f: \Y \to \RR$ belongs to the sub-gradient of a convex function if, and only if, $f$ is non-decreasing. When $f$ is continuous and $\Y$ is an open interval, a primitive $\psi$ of $f$ is indeed convex if, and only if, $\psi' = f$ is non-decreasing \cite[Proposition 17.7]{Bauschke:2017ki}. We now prove the result for more general $\Y$ and $f$. First, if $f(y) \in \partial \psi(y)$ for each $y \in \Y$ where $\psi: \RR \to \RR \cup \{+\infty\}$ is convex, then by \cite[Theorem 4.2.1 (i)]{Urruty96} $f$ is non-decreasing. To prove the converse  define $a:= \inf\{y: y \in \Y\}$, $I := (a,\infty)$ if $a \notin \Y$ (resp. $I := [a,\infty)$ if $a \in \Y$), and set $\bar{f}(x) := \sup_{y \in \Y, y \leq x} f(y) \in \RR \cup \{+\infty\}$ for each $x \in I$, $\bar{f}(x) = +\infty$, for $x \notin I$.  By construction $\bar{f}$ is non-decreasing. If $f$ is non-decreasing on $\Y$ then $\bar{f}(y) = f(y)$ for each $y \in \Y$ hence $\Y \subset \dom{\bar{f}} \subset I$ and $\dom{\bar{f}}$ is an interval. Choose an arbitrary  $b \in \Y$. As $\bar{f}$ is monotone it is integrable on each bounded interval one can define $\psi(x):= \int_{b}^{x} \bar{f}(t) dt$ for each $x  \in \dom{\bar{f}}$ (with the usual convention that if $x<b$ then $\int_{b}^{x} = -\int_{x}^{b}$) and $\psi(x) := +\infty$ for $x \notin \dom{\bar{f}}$.  Consider $x \in \dom{\bar{f}}$. Since $\bar{f}$ is non-increasing for $h \geq 0$ such that $x+h \in \dom{\bar{f}}$ we have $\psi(x+h)-\psi(x) = \int_{x}^{x+h} \bar{f}(t)dt \geq \bar{f}(x) h$; similarly for $h \leq 0 $ such that $x-h \in \dom{\bar{f}}$ we have $\psi(x)-\psi(x-h) = \int_{x-h}^{x} \bar{f}(t)dt \leq \bar{f}(x) (-h)$, hence $\psi(x-h)-\psi(x) \geq \bar{f}(x) h$. Combining both results shows $\psi(y)-\psi(x) \geq \bar{f}(x) (y-x)$ for each $x,y \in \dom{\bar{f}}$. This establishes that $\bar{f}(x) \in \partial \psi(x)$ for each $x \in \dom{\bar{f}}$, hence that $\psi$ is convex on its domain $\dom{\psi} = \dom{\bar{f}}$. To conclude, simply observe that for $y \in \Y \subset \dom{\bar{f}}$ we have $f(y) = \bar{f}(y) \in \partial \psi(y)$.
 \endproof

 \begin{example}[Quantization]\label{ex:QuantizationIsProx}
 In $\Y = [0,1) \subset \RR = \H$, consider $0 = x_{0} < x_{1} < \ldots < x_{q-1} < x_{q}=1$ and $v_{0} \leq \ldots \leq v_{q-1}$. Let $f$ be the quantization-like function so that $f(x) = v_{i}$ if and only if $x \in [x_{i},x_{i+1})$, for $0 \leq i < q$. Quantization traditionally corresponds to the case where $q \geq 2$ and for each $0 \leq i < q-1$, $x_{i+1}$ is the middle point between $v_{i}$ and $v_{i+1}$. Since $f$ is  non-decreasing, $f$ is the proximity operator of a function $\ph$. The image of $f$ is the discrete set of points $\{v_{0},\ldots,v_{q-1}\}$. 
\end{example}
Let us give another example to illustrate the role of the connectedness of the sets $\C$, $\C'$ in  \thm{ProxCh}.
\begin{example}\lab{ex:RoleConnectedness}
Consider the identity function $f(y):= y \mapsto y$ on a subset $\Y \subset \RR = \H$. 
Since $f$ is increasing it is a proximity operator by \cor{coscal}. Particular functions satisfying the equivalent properties \eq{tproa}, \eq{tprob} and 
\eq{tproc} of \thm{ProxCh} are $\ph_{0}: x \mapsto 0$,  $g_{0}: x \mapsto x^{2}/2$ and $\psi_{0}: y \mapsto y^{2}/2$. They further satisfy \eq{tggg} (resp. \eq{tpsi}) with $K=K'=0$ on $\RR$.
When $\Y \subset \RR$ is polygonally connected, $\Ima{f}$ is also polygonally connected by the continuity of $f$ and \thm{ProxCh} implies that $\ph_{0},g_{0},\psi_{0}$ are, up to global additive constants $K,K'$, the only functions satisfying \eq{tproa}, \eq{tprob} and 
\eq{tproc}.
Now, consider as a particular example of disconnected set $\Y = (-\infty,0) \cup (1,+\infty)$.
We exhibit two other functions $g,\psi$ such that $\ph_{0},g,\psi$ also satisfy \eq{tproa}, \eq{tprob} and 
\eq{tproc}, but  \eq{tggg} fails on the disconnected set $\C := \Ima{f} = \Y$ (resp. \eq{tpsi} fails on the disconnected set $\C':=\Y$).
Intuitively, what happens is that the presence of a ``hole'' (the interval $[0,1]$) in $\Y$ gives some freedom in designing separately the components of these functions on each connected component.
For this, consider $H: [0,1] \to [0,1]$ any continuous increasing function such that $H(0) = 0$, $H(1) = 1$ and $C:=\int_0^{1} H(t)dt \neq 1/2$. Observe that the function
\[
h(x) := \begin{cases}
\tfrac{x^{2}}{2}, x < 0\\
\int_{0}^{x} H(t) dt, 0 \leq x \leq 1\\
\int_{0}^{1} H(t) dt + \tfrac{x^{2}-1}{2}, x > 1.
\end{cases}.
\]
is convex and satisfies $\partial h(x) = \{h'(x)\} = \{x\}$ for each $x \in \Y$. As a result the functions $g := h$ and $\psi := h$ also satisfy properties \eq{tprob} and 
\eq{tproc} of \thm{ProxCh}. Yet on the interval $(-\infty,0)$ we have $g(x) = g_{0}(x) = \ph_{0}(x) + \tfrac{x^{2}}{2} + K_{0}$ with $K_{0}=0$,  while on the interval $(1,+\infty)$ we have $g(x) = g_{0}(x)+C-1/2 = \ph_{0}(x) + \tfrac{x^{2}}{2} + K_{1}$ with $K_{1} = C-1/2 \neq 0 = K_{0}$. Similarly $\psi(x)-\psi_{0}(x)$ is not constant on $\Y$. This shows that \eq{tggg} (resp. \eq{tpsi}) fails to hold on $\C := \Ima{f}$ (resp. $\C' := \Y$). 
\end{example}

Consider now functions $f: \RR^{n} \to \RR^{n}$ given by $f(y) = (f_{i}(y))_{i=1}^{n}$. 
When each $f_{i}$ can be written as $f_{i}(y) = h_{i}(y_{i})$, the function is said to be separable. If each $h_{i}$ is a scalar proximity operator then the function $f$ is also a proximity operator, and vice-versa. This can be seen, e.g., by writing $h_{i} = \prox_{\ph_{i}}$ and $f = \prox_{\ph}$ with $\ph(x) := \sum_{i=1}^{n}\ph_{i}(x_{i})$. 
All examples below hold for the components of separable functions.

As recalled in~\prop{MoreauNE} it is known \cite[Proposition 2.4]{Combettes:2007aa} that a function $f: \RR \to \RR$ is the proximity operator of a {\em convex l.s.c.} penalty $\ph$ if, and only if, $f$ is nondecreasing and nonexpansive: $|f(y)-f(y')| \leq |y-y'|$ for each $y,y' \in \RR$. 

A  particular example is that of scalar {\em thresholding rules} which are known \cite[Proposition 3.2]{Antoniadis:2007aa} to be the proximity operator of a ({\em continuous positive}) penalty function. As we will see in Section~\ref{sec:social},~\thm{ProCon} also allows to characterize whether certain {\em block-thresholding rules} \cite{Hall:1997wz,Cai:2001tj,Kowalski:2013vx} are proximity operators.
 
Our next example illustrates the functions appearing in~\thm{ProCon} on the classical hard-thresholding operator, which is the proximity operator of a nonconvex function.
\begin{example}[Hard-thresholding]\label{ex:L0}
In $\Y=\H = \RR$ consider $\lambda>0$ and the weighted $\ell^{0}$ penalty
\[
\ph(x) := \begin{cases}0,\ \text{if}\ x=0;\\ \lambda,\ \text{otherwise}.\end{cases}
\]
Its (set-valued) proximity operator is
\[
\prox_\ph(y)=\lp\{\barr{lll} \{0\} &  \mb{if} & |y| < \sqrt{2\la}\\
\{0,\sqrt{2\la} \} & \mb{if} &\ y = \sqrt{2\la}\\
\{-\sqrt{2\la} ,0\} &\mb{if}& y = -\sqrt{2\la}\\
\{y\}& \mb{if} & |y|> \sqrt{2\la}\earr \rp.
\]
which is discontinuous. Choosing $\pm \sqrt{2\lambda}$ as the value at $y = \pm \sqrt{2\lambda}$ yields a function $f(y) \in \prox_\ph(y)$ with disconnected (hence nonconvex) range $\qu \Ima{f} = (-\infty,-\sqrt{2\la}] \cup \{0\} \cup [\sqrt{2\la},+\infty)$,
\[
f(y) := \begin{cases}0,\ \text{if}\ |y| < \sqrt{2\lambda}\\ y,\ \text{if}\ |y| \geq \sqrt{2\lambda}\end{cases}
\qquad
\]
Since $\Y$ is convex, the potential $\psi$ is characterized by \eq{psivsphi}. For $K:=0$ we get
\[
\psi(y) = yf(y)-\tfrac{1}{2} f^2(y)-\ph(f(y))
=
\begin{cases} 0,\ \text{if}\ |y| < \sqrt{2\lambda}\\ 
y^2/2-\lambda,\ \text{otherwise}\end{cases}
=\max(y^2/2-\lambda,0).
\]
This is indeed a convex potential, and $f(y) \in \d\psi(y)$ for each $y \in \RR$. 
\end{example}

Our last example of this section is a scaled version of soft-thresholding: it is still a proximity operator, however for $C>1$ the corresponding penalty is nonconvex, and is even unbounded from below.
\begin{example}[Scaled soft-thresholding]\label{ex:GCvxPhiNot}
In $\Y=\H = \RR$ consider 
\[
f(y) := 
\begin{cases}
0,\ \text{if}\ |y| < 1\\
C(y-1),\ \text{if}\ y\geq 1\\
C(y+1),\ \text{if}\ y\leq -1
\end{cases}
= C y \max(1-1/|y|,0).
\]
This function has the same shape as the classical soft-thresholding operator, but is scaled by a multiplicative factor $C$. 
When $C=1$, $f$ is the soft-thresholding operator which is the proximity operator of the absolute value, $\ph(x) = |x|$, which is convex. For $C>1$, as $f$ is expansive, by \prop{MoreauNE} it cannot be the proximity operator of any convex function. Yet, as $f$ is monotonically increasing, $f(y)$ is a subgradient of its ``primitive'' $\psi(y) = \tfrac{C}{2} \lp(\max(|y|-1,0)\rp)^2 = \tfrac{C}{2} y^2 \lp(\max(1-1/|y|,0)\rp)^2 = \tfrac{f^2(y)}{2C}$ which is convex. Moreover, by \cor{coscal}, $f$ is still the proximity operator of some (necessarily nonconvex) function $\ph(x)$. By \eq{psivsphi}, up to an additive constant $K \in \RR$, $\ph$ satisfies 
\[
\ph(f(y)) =  yf(y)-\tfrac{1}{2} f^2(y)-\psi(y)
= yf(y) - \tfrac{1+C}{2C}f^2(y), \all y \in \RR
\]
For $x > 0$, writing $x = f(y)$ with $y = f^{-1}(x) = 1+x/C$ yields $\ph(x) = \ph(f(y)) = (1+x/C)x-\tfrac{1+C}{2C}x^2$.
Similar considerations for $x<0$ and for $x=0$ show that $\ph(x) = |x|+\left(\tfrac{1}{C}-1\right)\tfrac{x^{2}}{2}$. When $C>1$, $\ph$ is indeed not bounded from below, and not convex.
\end{example}

\section{When is social shrinkage a proximity operator ?} \lab{sec:social}
We conclude this paper by studying so-called social shrinkage operators, which have been introduced to mimic classical sparsity promoting proximity operators when certain types of structured sparsity are targeted. We show that the characterization of proximity operators obtained in this paper provides answers to questions raised by Kowalski et al \cite{Kowalski:2013vx} and by Varoqueaux et al  \cite{Varoquaux:2016wd} on the link between such non-separable shrinkage operators and proximity operators.

Most proximity operators are indeed not separable.
A classical example is the proximity operator associated to mixed $\ell_{12}$ norms, which enforces group-sparsity.
\begin{example}[Group-sparsity shrinkage] \lab{ex:GroupSparsity}
Consider a partition $\mathcal{G} = \{G_{1},\ldots,G_{p}\}$ of $\set{n}$, the interval of integers from $1$ to $n$, into disjoint sets called {\em groups}. 
Let $x_{G}$ be the restriction of $x \in \RR^{n}$ to its entries indexed by $G \in \mathcal{G}$, and define the {\em group $\ell^{1}$ norm}, or {\em mixed $\ell_{12}$ norm}, as 
\beq \lab{GL}
\ph(x) := \sum_{G \in \mathcal{G}} \|x_G\|_{2}.
\eeq
The proximity operator $f(y) := \prox_{\lambda \ph}$ is the {\em group-sparsity shrinkage} operator with threshold $\lambda$
\beq \lab{GLprox}
\forall i \in G,\quad f_{i}(y):= y_{i} \left(1-\frac{\lambda}{\|y_{G}\|_{2}}\right)_{+}.
\eeq
\end{example}
The group-LASSO penalty~\eq{GL} appeared in statistics in the thesis of Bakin \cite[Chapter 2]{Bakin1999}. It was popularized by Yuan and Lin \cite{Yuan:2006jy} who introduced an iterative shrinkage algorithm to address the corresponding optimization problem. A generalization is Group Empirical Wiener / Group Non-negative Garrotte, see e.g. \cite{Fevotte:2015fb}
\beq \lab{GEWprox}
\forall i \in G,\quad f_{i}(y):= y_{i} \left(1-\frac{\lambda^2}{\|y_{G}\|_{2}^2}\right)_{+},
\eeq
see also  \cite{Antoniadis:2007aa} for a review of thresholding rules, and \cite{bach:hal-00613125} for a review on sparsity-inducing penalties.

To account for varied types of structured sparsity, \cite{KOWALSKI-2009-315855,Kowalski:2009vf} empirically introduced the so-called Windowed Group-LASSO. A weighted version for audio applications was further developed in \cite{Siedenburg:2011uq} which coins the notion of {\em persistency}, and the term {\em social sparsity} was coined in \cite{Kowalski:2013vx} to cover Windowed Group-LASSO, as well as other structured shrinkage operators. As further described in these papers, the main motivation of such social shrinkage operators is to obtain flexible ways of taking into account (possibly overlapping) {\em neighborhoods} of a coefficient index $i$ rather than {\em disjoint groups} of indices to decide whether or not to set a coefficient to zero. These are summarized in the definition below.

\BD[Social shrinkage]\lab{defsoc}
Consider a family $N_{i} \subset \set{n}$, $i \in \set{n}$ of sets such that $i \in N_{i}$. The set $N_i$ is called a {\em neighborhood} of its index $i$. Consider nonnegative weight vectors $w^{i}= (w^{i}_{\ell})_{\ell=1}^{n}$ such that $\supp(w^i) = N_i$. Windowed Group Lasso (WG-LASSO) shrinkage is defined as $f(y) := (f_{i}(y))_{i=1}^{n}$ with 
\begin{equation}\lab{WG-LASSO}
\forall i,\quad f_{i}(y) := y_{i} \left(1-\frac{\lambda}{\|\diag(w^{i}) y\|_{2}}\right)_{+}
\end{equation}
and Persistent Empirical Wiener (PEW) shrinkage (see \cite{Siedenburg:2014ih} for the unweighted version) with
\begin{equation}\lab{PEW}
\forall i,\quad f_{i}(y) := y_{i} \left(1-\frac{\lambda^{2}}{\|\diag(w^{i}) y\|_{2}^{2}}\right)_{+}.
\end{equation}
\ED
Kowalski et al \cite{Kowalski:2013vx} write ``{\em while the classical proximity operators\footnote{that are explicitly constructed as the proximity operator of a convex l.s.c. penalty, e.g., soft-thresholding.} are directly linked to convex regression problems with mixed norm priors on the coefficients, 
[those]
 new, structured, shrinkage operators cannot be directly linked to a convex minimization problem}''. Similarly,  Varoquaux et al \cite{Varoquaux:2016wd} write that Windowed Group Lasso ``{\em is not the proximal operator of a known penalty}''. They leave open the question of whether social shrinkage is the proximity operator of some yet to be discovered penalty. Using~\thm{jac1}, we answer these questions for generalized social shrinkage operators.
The answer is negative unless the involved neighborhoods form a partition.

\BD[Generalized social shrinkage] \lab{gss}
Consider subsets $N_{i} \subset  \set{n}$ and nonnegative weight vectors $w^{i} \in \RR_{+}^{n}$ such that $i \in N_{i}$ and $\supp(w^{i}) = N_{i}$ for each $i \in \set{n}$.
Consider $\lambda > 0$ and a family of $C^1(\RR_+^*)$ scalar functions $h_{i}$, $i \in \set{n}$ such that $h'_{i}(t) \neq 0$ for $t \in \RR_+^*$. 
A generalized social shrinkage operator 
is defined as $f(y) := (f_{i}(y))_{i=1}^{n}$ with 
\[ 
f_{i}(y) :=  
\begin{cases} y_{i} h_{i}
\left(\|\diag(w^{i})y\|_2^{2}\right),\ \text{if}\ \|\diag(w^{i})y\|_2> \lambda,\\ 0 \ \text{otherwise}. \end{cases} 
\] 
\ED
We let the reader check that the above definition covers Group LASSO \eq{GLprox}, Windowed Group-LASSO~\eq{WG-LASSO}, Group Empirical Wiener \eq{GEWprox} and Persistent Empirical Wiener shrinkage~\eq{PEW}.

\BL\label{le:PEWNotProx}
Let $f: \RR^n \to \RR^n$ be a generalized social shrinkage operator and $N_i \subset \set{n}$, $w^i \in \RR^n_+$, $i \in \set{n}$ be the corresponding families of neighborhoods and weight vectors. 
If $f$ is a proximity operator then there exists a partition $\mathcal{G} = \{G_{p}\}_{p=1}^{P}$ of the set $\set{n}$ of indices such that: for each $p$ and all $i,j \in G_{p}$ we have $w^{i}=w^{j}$ and $\supp(w^{i}) = \supp(w^{j}) = G_{p}$. As a consequence for $i \in G_{p}$, $j \in G_{q}$ with $p \neq q$, the weight vectors $w^{i}$ and $w^{j}$ have disjoint support.
\EL

The proof of \lem{le:PEWNotProx} is postponed to Appendix~\ref{app:PfLePEWNotProx}.
An immediate consequence of this lemma is that if $f$ is a generalized social shrinkage operator, then the neighborhood system $N_{i} = \supp(w^{i})$ coincides with the groups $G$ from the partition $\mathcal{G}$. In particular, the neighborhood system must form a partition. By contraposition we get the following corollary:
\BC\lab{cor:PEWNotProx}
Consider non-negative weights $\{w^i\}$ as in \defi{gss} and $\{N_i\}$ the corresponding neighborhood system.
Assume that there exists $i,j$ such that $N_i \neq N_j$ and $N_i \cap N_j \neq \void$.
\bit
\item Let $f$ be the WG-LASSO shrinkage \eq{WG-LASSO}. There is no penalty $\ph$ such that $f = \prox_\ph$.
\item Let $f$ be the PEW shrinkage \eq{PEW}. There is no penalty $\ph$ such that $f = \prox_\ph$.
\eit
\EC
In other words, WG-LASSO / PEW can be a proximity operator {\em only if} the neighborhood system has {\em no overlap}, i.e. with ``plain'' Group-LASSO \eq{GLprox} /  Group Empirical Wiener \eq{GEWprox}.

\section*{Acknowledgements}
The first author wishes to thank Laurent Condat, Jean-Christophe Pesquet and Patrick-Louis Combettes for their feedback that helped improve an early version of this paper, as well as the anonymous reviewers for many insightful comments that improved it much further.

\appendix

\section{Proofs} 

The proofs of technical results of Section~\ref{sec:prox} are provided in Sections~\ref{pf:GProxCh} (\thm{GProxCh}), \ref{lsp} (\lem{SdCont}), \ref{pfuniqueC1} (\cor{uniqueC1}) and \ref{app:PfLePEWNotProx} (\lem{le:PEWNotProx}). As a preliminary we give brief reminders on some useful but classical notions in Sections~\ref{frechet}-\ref{sec:subdiffunique}.

\subsection{Brief reminders on (Fr{\'e}chet) differentials and gradients in Hilbert spaces}\lab{frechet}
Consider $\H,\H'$ two Hilbert spaces.
A function $\ff: \X \to \H'$ where $\X \subset \H$ is an open domain is (Fr{\'e}chet) differentiable at $x$ if there exists a continuous linear operator $L:\H\to\H'$ such that $\lim_{h \to 0}\|\ff(x+h)-\ff(x)-L(h)\|_{\H'}/\|h\|_{\H} = 0$. The linear operator $L$ is called the differential  of $\ff$ at $x$ and denoted $D\ff(x)$.  When $\H' = \RR$, $L$ belongs to the dual of $\H$, hence there is $u \in \H$ --called the gradient of $\ff$ at $x$ and denoted $\nabla \ff(x)$-- such that $L(h) = \<u,h\>,\ \all h \in \H$.

\subsection{Subgradients and subdifferentials for possibly nonconvex functions} \lab{sec:lemmasProx}
We adopt a gentle definition which is familiar when $\ff$ is a convex function. Although this is possibly less well-known by non-experts, this definition is also valid when $\ff$ is possibly nonconvex, see e.g. \cite[Definition 16.1]{Bauschke:2017ki}.
\BD\lab{nlsou} Let $\ff:\H\to \RR\cup \{+\infty\}$ be a proper function. 
The {\em \sd} $\bd \ff(x)$ of $\ff$ at $x$ is the set of all 
$u\in\H$, called {\em \sg}s of $\ff$ at $x$, such that 
\beq\label{DSG}
\ff(x') \geq \ff(x) + \< u,x'-x\>,\quad \forall x' \in \H.
\eeq
If $x\not \in \dom{\ff}$, then  $\bd \ff(x)=\void$.
The function $\ff$ is {\em \sda} at $x \in \H$ if $\bd\ff(x) \neq \void$. 
The domain of $\bd\ff$ is $\dom{\bd\ff} := \{x \in \H, \bd\ff(x) \neq \void\}$. It satisfies $\dom{\bd\ff} \subset \dom{\ff}$.
\ED
\BF\lab{finf} When $\bd\ff(x)\neq\void$  
the inequality in \eq{DSG} is trivial for each $x'\not\in\dom{\ff}$ \db since it amounts to $+\infty = \ff(x')-\ff(x) \geq \< u,x'-x\>$.
\EF

 \defi{nlsou} leads to
the well-known Fermat's rule \cite[Theorem 16.3]{Bauschke:2017ki}

\BT\lab{fe} Let $\ff:\H\to \RR\cup \{+\infty\}$ be a proper function.
A point $x\in \dom{\ff}$ is a global minimizer of $\ff$ if and only if
\[ 
0 \in \bd \ff(x).
\]
\ET

If $\ff$ has a global minimizer at $x$, then by  \thm{fe} the set $\bd \ff(x)$ is non-empty. However, $\bd \ff(x)$ can be empty, e.g., at local minimizers that are not the global minimizer:  
\BE\lab{sin} Let $\ff(x)=\frac12x^2-\cos(\pi x)$. 
The global minimum of $\ff$ is reached at $x=0$ where $\bd \ff(x)= f'(x)=0$.
At $x=\pm 1.7 \bar{9} $ $\ff$ has local minimizers where $\bd \ff(x)=\void$ (even though $\ff$ is $\C^\infty$).
For $|x|<0.53$ one has $\bd\ff(x)=\nabla\ff(x)$ with $\ff''(x)\geq0$ and  for $0.54 < |x| < 1.91$ $\bd\ff(x)=\void$. 
\EE
The proof of the following lemma is a standard exercice in convex analysis \cite[Exercice 16.8]{Bauschke:2017ki}.
\BL\lab{cac} Let $\ff:\H\to \RR\cup \{+\infty\}$ be a proper function such that (a) $\dom{\ff}$ is convex
 and (b)
$\bd \ff(x)\neq \void$ for each  $x \in \dom{\ff}$.
Then $\ff$ is a convex function. 
\EL

\BD (Lower convex envelope of a function)\lab{cvxlow}\\
Let  $\ff: \H \to \RR \cup \{+\infty\}$ be proper with $\dom{\bd\ff} \neq \void$. Its lower convex envelope,\footnote{also known as convex hull, \cite[p. 57]{Rockafellar:1998ge},\cite[Definition 2.5.3]{Urruty96} }
denoted $\breve{\ff}$, is the pointwise supremum of all the convex lower-semicontinuous functions minorizing $\ff$
\beq
\lab{deflcev}
\breve{\ff}(x) := \sup \{\gg(x)\, |\, \gg: \H \to \RR \cup \{+\infty\}, \gg\ \mb{convex l.s.c.},\ \gg(z) \leq \ff(z), \all z \in \H \},\qu \all x \in \H.
\eeq
The function $\breve{\ff}$ is proper, convex and lower-semicontinuous. 
It satisfies 
\beq\lab{cvxenvin}
\breve{\ff}(x) \leq \ff(x), \all x \in \H.
\eeq
\ED

\BP\lab{sdcvxenv}
Let $\ff: \H \to \RR \cup \{+\infty\}$ be proper with $\dom{\bd\ff} \neq \void$.
For any $x_{0} \in \dom{\bd\ff}$ we have
 $\breve{\ff}(x_0) = \ff(x_0)$, $\bd\ff(x_0) = \d\breve{\ff}(x_0)$.
\EP
\proof
As $\bd \ff(x_0) \neq \void$, by \cite[Proposition 13.45]{Bauschke:2017ki}, $\breve{\ff}$ is the so-called biconjugate $\ff^{**}$ of $\ff$ \cite[Definition 13.1]{Bauschke:2017ki}. Moreover, \cite[Proposition 16.5]{Bauschke:2017ki} yields $\ff^{**}(x_{0}) = \ff(x_{0})$ and $\d\ff^{**}(x_{0}) = \d\ff(x_{0})$.
\endproof

We need to adapt \cite[Proposition 17.31]{Bauschke:2017ki} to the case where $\ff$ is proper but possibly nonconvex, with a stronger assumption of Fr{\'e}chet (instead of G{\^a}teaux) differentiability.
\BP\lab{sdisg}
If $\bd\ff(x) \neq \void$ and $\ff$ is (Fr\'echet) differentiable at $x$ then $\bd\ff(x) = \{\nabla \ff(x)\}$.
\EP

\proof
Consider $u \in \bd\ff(x)$. As $\ff$ is differentiable at $x$ there is an open ball $\B$ centered at $0$ such that $x+h \in \dom{\ff}$ for each $h \in \B$. For each $h \in \B$,  \defi{nlsou} yields
 \[
 \ff(x-h)-\ff(x) \geq \<u,-h\>\qu\qu \mb{and}\qu\qu
 \ff(x+h)-\ff(x) \geq \<u,h\>
 \]
 hence $-(\ff(x-h)-\ff(x)) \leq \<u,h\> \leq \ff(x+h)-\ff(x)$. Since $\ff$ is Fr{\'e}chet differentiable at $x$, letting $\|h\|$ tend to zero yields 
 \[
 -\left(\<\nabla \ff(x),-h\> + o(\|h\|)\right) \leq \<u,h\> \leq \<\nabla \ff(x),h\> + o(\|h\|)
 \]
 hence $\<u-\nabla\ff(x),h\> = o(\|h\|)$, $\all h\in\B$. This shows that $u = \nabla \ff(x)$.
\endproof

\subsection{Characterizing functions with a given subdifferential}\lab{sec:subdiffunique}

\cor{uni} below generalizes a result of Moreau \cite[Proposition 8.b]{Moreau65} characterizing functions by their subdifferential. It shows that one only needs the subdifferentials to intersect. We begin in dimension one.

\BL\lab{uniscal} 
Consider $a_0,a_1: \RR \to \RR \cup \{+\infty\}$ convex functions such that $\dom{a_i} = \dom{\bd a_i} = [0,1]$ and $\d a_0(t) \cap \d a_1(t) \neq \void$ on $[0,1]$. 
Then there exists a constant $K \in \RR$ such that $a_1(t)-a_0(t)=K$ on $[0,1]$. 
\EL
\proof 
As $a_i$ is convex it is continuous on $(0,1)$ \cite[Theorem 3.1.1, p16]{Urruty96}. Moreover, by \cite[Proposition 3.1.2]{Urruty96} we have $a_{i}(0) \geq \lim_{t \to 0, t>0} a_{i}(t) =:a_{i}(0_+)$, and since  $\bd a_{i}(0) \neq \void$, there is $u_{i} \in \bd a_{i}(0)$ such that $a_{i}(t) \geq a_{i}(0) + u_{i}(t-0)$ for each $t \in [0,1]$ hence $a_{i}(0_+) \geq a_{i}(0)$. This shows that $a_{i}(0_+) = a_{i}(0)$, and similarly $ \lim_{t \to 1, t<1} a_{i}(t) = a_{i}(1)$, hence $a_{i}$ is continuous on $[0,1]$ relatively to $[0,1]$. In addition, $a_{i}$ is differentiable on $[0,1]$ except on a countable set $B_i \subset [0,1]$ \cite[Theorem 4.2.1 (ii)]{Urruty96}.

For $t \in [0,1] \backslash (B_{0} \cup B_{1})$ and $i \in \{0,1\}$, \prop{sdisg} yields $\d a_i(t) = \{a'_i(t)\}$ hence the function $\delta := a_1-a_0$ is continuous on $[0,1]$ and differentiable on $[0,1] \backslash (B_0 \cup B_1)$. 
For $t \in I \backslash (B_0 \cup B_1)$, $\{a'_0(t)\} \cap \{a'_1(t)\} = \d a_0(t) \cap \d a_1(t) \neq \void$, hence $a'_0(t)=a'_1(t)$ and $\delta'(t) = 0$. A classical exercice\footnote{for a proof see e.g. (in french) \url{https://fr.wikipedia.org/wiki/Lemme\_de\_Cousin} section 4.9, version from 13/01/2019.} in real analysis \cite[Example 4]{Thomson:2007ui} is to show that if a function $f$ is continuous on an interval, and differentiable with zero derivative except on a countable set, then $f$ is constant.
As $B_0 \cup B_1$ is countable it follows $\delta$ is constant on $(0,1)$. As it is continuous on $[0,1]$, it is constant on $[0,1]$.
\endproof

\BC \lab{uni}
Let $\ff_0,\ff_1: \H \to \RR \cup \{+\infty\}$ be proper and $\C \subset \H$ a non-empty polygonally connected set.
Assume that for each $z \in \C$, $\bd \ff_0(z) \cap \bd \ff_1(z) \neq \void$; then there is a constant $K \in \RR$ such that $\ff_1(x) -\ff_0(x) = K$, $\all x \in \C$.
\EC
\BR
Note that the functions $\ff_{i}$ and the set $\C$ are not assumed to be convex.
\ER
\proof 
The proof is in two parts.

(i)
Assume that $\C$ is convex and fix some $x^* \in \C$. 
Consider  $x \in \C$, and define $a_i(t):= \ff_i(x^*+t(x-x^*))$, for $i=0,1$ and each $t \in [0,1]$, and $a_i(t)=+\infty$ if $t\not\in[0,1]$. 
As $\C$ is convex, $z_t := x^*+t(x-x^*) \in \C$ hence 
for each $t \in [0,1]$ there exists $u_t \in \bd \ff_0(z_t) \cap \bd\ff_1(z_t)$. 
By \defi{nlsou} for each $t,t' \in [0,1]$,
\[
a_i(t')-a_i(t) 
= \ff_i(x^*+t'(x-x^*))-\ff_i(x^*+t(x-x^*)) 
\geq \<u_t,(t'-t)(x-x^*)\>
=\<u_t,x-x^*\>(t'-t).
\]
For $t \in [0,1]$ and $t' \in \RR \backslash [0,1]$ since $a_{i}(t') = +\infty$ the inequality $a_i(t')-a_i(t) \geq \<u_t,x-x^*\>(t'-t)$ also obviously holds, hence $\<u_t,x-x^*\> \in \bd a_i(t)$, $i=0,1$. Thus $\bd a_i(t)\neq\void$ for each $t \in [0,1]$, so by \lem{cac}  $a_i$ is convex on $[0,1]$ for $i=0,1$, and $\<u_t,x-x^*\> \in \d a_0(t) \cap \d a_1(t)$ for each $t \in [0,1]$. By \lem{uniscal}, there exists $K \in \RR$ such that $a_1(t)-a_0(t)= K$ for each $t \in [0,1]$. Therefore, 
\[
\ff_1(x)-\ff_0(x) = a_1(1)-a_0(1) = a_1(0)-a_0(0) = \ff_1(x^*)-\ff_0(x^*) = K.
\]
As this holds for each $x \in \C$, we have established the result as soon as $\C$ is convex. 

(ii) Now we prove the result when $\C$ is polygonally connected. Fix some $x^* \in \C$ and define $K:= \ff_1(x^*)-\ff_0(x^*)$. Consider $x \in \C$: by the definition of polygonal connectedness, there exists an integer $n \geq 1$ and $x_j \in \C$, $0 \leq j \leq n$ with $x_0 = x^*$ and $x_n = x$ such that the (convex) segments $\C_j = [x_j,x_{j+1}] = \{t x_j + (1-t) x_{j+1}, t \in [0,1]\}$ satisfy $\C_j \subset \C$. Since each $\C_j$ is convex, the result  established in (i) implies that $\ff_1(x_{j+1})-\ff_0(x_{j+1}) = \ff_1(x_j)-\ff_0(x_j)$  for $0 \leq j < n$. 
This shows that $\ff_1(x)-\ff_0(x) = \ff_1(x_{n})-\ff_0(x_{n}) = \ldots = \ff_1(x_{0})-\ff_0(x_{0}) = \ff_1(x^*)-\ff_0(x^*) =K$. 
\endproof

\subsection{Proof of \thm{GProxCh}} \lab{pf:GProxCh}

The indicator function of a set $\S$ is denoted
\[ \chi_{\S}(x):=\lp\{ \barr{lll} 0&\mb{if}&x\in \S, \\ +\infty&\mb{if}& x\not\in \S.\earr\right.\]

\paragraph{\eq{proa} $\Rightarrow$ \eq{prob}} 
We introduce the function $\ff:\H \to\RRinf$ by  
\beq
\lab{gggtmp}  \ff := b+\ph + \chi_{\Ima{f}}. 
\eeq
Consider $x \in \Ima{f}$. By definition $x= f(y)$ where $y \in \Y$, hence by \eq{proa} $x$ is a global minimizer of $x'\mapsto \lp\{D(x',y)+\ph(x')\rp\}$.
Therefore, we have
 \beq
 \all x' \in \H,\ -\<A(y),x'\>+\underbrace{b(x')+\ph(x') + \chi_{  \Ima{f}  } (x')}_{=\ff(x')}
 \geq 
  -\<A(y),x\>+\underbrace{b(x)+\ph(x) + \chi_{\Ima{f}}(x)}_{=\ff(x)} \lab{amphi}
 \eeq
which is equivalent to
\beqn
\all x'\in \H\qu \qu \ff(x') \geq \ff(x) + \<A(y),x'-x\>    \lab{gsd}
\eeqn
meaning that  
$A(y)\in \bd \ff\lp(f(y)\rp)$.
As this holds for each $y \in \Y$ such that $f(y)=x$, we get $A(f^{-1}(x))\subset \bd \ff(x)$. 
Consider $g_1 := \breve \ff $ according to \defi{cvxlow}.
Since $g_1$ is convex l.s.c. and 
\beqn
\all x \in \Ima{f}, \bd\ff(x) \neq \void, \lab{bffsd}
\eeqn
by \prop{sdcvxenv}, $\bd \ff(x) = \d g_1(x)$ and $\ff(x) = g_{1}(x)$ for each $x \in \Ima{f}$. This establishes \eq{prob} with $g := g_1 = \breve{\ff}$.\\

\paragraph{\eq{prob} $\Rightarrow$ \eq{proa}}
Set $\ff_1: = g+\chi_{\Ima{f}}$. 
By \eq{prob}, $\d g(x) \neq \void$ for each $x \in \Ima{f}$.
Since $\dom{\d g} \subset \dom{g}$ it follows that $\Ima{f} \subset \dom{g}$ and consequently 
\[
\dom{\ff_1}=\Ima{f}.
\]
Consider $y \in \Y$ and $x:=f(y)$ so that $x \in \Ima{f}$, hence $\ff_1(x)=g(x)$ and $A(y)\in A(f^{-1}(x)) \subset \d g(x)$ 
where the inclusion comes from \eq{prob}. It follows that for each $(x,x')\in \Ima{f} \times \H$ one has 
\[
\ff_1(x')= g(x')+\chi_{\Ima{f}}(x') \geq g(x') \geq g(x)+\<A(y),x'-x\> = \ff_{1}(x)+\<A(y),x'-x\>,
\]
showing that $A(y)\in \bd \ff_1(x)$. This is equivalent to~\eq{gsd} with $\ff := \ff_1$, and since $\dom{\ff_1} = \Ima{f}$, the inequality in~\eq{amphi} holds with $\ph(x) := \ff_1(x)-b(x)$, i.e., $x$ is a global minimizer of $D(x',y)+\ph(x')$. Since this holds for each $y \in \Y$, this establishes~\eq{proa} with $\ph := \ff_1-b = g-b+\chi_{\Ima{f}}$.\\

 \paragraph{\eq{pro2}} 
Consider $\ph$ and $g$ satisfying \eq{proa} and \eq{prob},  respectively. 
Let\footnote{In general, we may have $g\neq g_1$ as there is no connectedness assumption on $\dom{\ff}$.} $g_1 := \breve{\ff}$ with $\ff$ defined in \eq{gggtmp}. 
Following the arguments of \eq{proa} $\Rightarrow$ \eq{prob} we obtain that $g_1$ (just as $g$) satisfies \eq{prob}. For each $x \in \C$ we thus have $\d g(x) \cap \d g_1(x) \supset A(f^{-1}(x)) \neq \void$ with $g,g_1$ convex l.s.c. functions. Hence, by \cor{uni}, since $\C$ is polygonally connected, there is a constant $K$ such that $g(x) = g_1(x)+K$, $\all x \in \C$. 
To establish the relation \eq{ggg} between $g$ and $\ph$ we now show that $g_1(x) = b(x) + \ph(x)$ on $\C$. By \eq{bffsd} and \prop{sdcvxenv} we have
$\breve\ff(x)=\ff(x)$ for each $x \in \Ima{f}$, hence as $\C \subset \Ima{f}$ we obtain $g_1(x) := \breve{\ff}(x) = \ff(x) = b(x)+\ph(x)$ for each $x \in \C$. This establishes \eq{ggg}. \\

\paragraph {\eq{proabis} $\Rightarrow$ \eq{proc}} 
Define
\beqn
\gg(y) &:=&
\begin{cases}
+\infty, &\forall y \notin \Y\\
\<B(f(y)),y\> -b(f(y))- \ph(f(y)),\quad&\forall y \in \Y.
\end{cases} \lab{psitmp}
\eeqn
Consider $y \in \Y$. From \eq{proabis}, for each $y'$ the global minimizer of $x \mapsto {\wt{D}}(x,y')+\ph(x)$ is reached at $x'=f(y')$. 
Hence, for $x = f(y)$ we have
\[ -\<B(f(y')),y'\>+b(f(y'))+\ph(f(y'))  \leq 
-\<B(x),y'\>+b(x)+\ph(x)
=
-\<B(f(y)),y'\>+b(f(y))+\ph(f(y))  \]
Using this inequality we obtain that
\beqnn \all  y'\in\Y,\
\gg(y')-\gg(y) &=& 
-\<B(f(y)),y\> + b(f(y)) + \ph(f(y)) + \<B(f(y')),y'\>-b(f(y'))-\ph(f(y')) \\
& \geq &
\<B(f(y)),y'\>-\<B(f(y)),y\> \geq \<  B(f(y)), y'-y\>
\eeqnn 
This shows that 
\beqn
B(f(y)) \in \bd\gg(y).\lab{procrhosd}
\eeqn
Set $\psi_1 := \breve \gg $ according to \defi{cvxlow}. 
Then the function  $\psi_1$ is convex l.s.c. and for each $y \in \Y$ the function $ B(f(y))$ is well defined, 
so $\bd\gg(y) \neq \void$. Hence, by \prop{sdcvxenv}, $\bd\gg(y) =\d\breve\gg(y)= \d\psi_1(y)$ and $\gg(y)=\breve{\gg}(y) = \psi_{1}(y)$ for each $y \in \Y$. This establishes \eq{proc} with $\psi := \psi_1 =\breve{\gg}$.\\

\paragraph{\eq{proc} $\Rightarrow$ \eq{proabis}} 
Define $h:\Y\to \RR$ by
\[  
h(y):=  \< B(f(y)),y\> - \psi(y)
\]
Since $B(f(y')) \in \d\psi(y')$ 
with $\psi$ convex by \eq{proc}, applying \defi{nlsou} to $\d \psi$ 
yields $\psi(y) - \psi(y') \geq \<y-y',B(f(y'))\>$. 
Using this inequality, one has
\beq \barr{lll} \all y,y'\in\Y\qu\qu  h(y')-h(y) &=& \<B(f(y')),y'\> - \psi(y') -\<B(f(y)),y\>+ \psi(y)  \\&\geq&
\<B(f(y')),y'\> - \<B(f(y)),y\> + \<B(f(y')),y-y'\>\\
&=&  \big\<B(f(y') -B(f(y)),\ y\big\>  
\earr \lab{22}
\eeq
Noticing that for each $x \in \Ima{f}$ there is $y \in \Y$ such that $x=f(y)$, we can define $\ff:\H\to \RRinf$ obeying $\dom{\ff}=\Ima{f}$ by
\[ 
\ff(x) := \lp\{\barr{lll} h(y ) &\mb{ with}~~ y\in f^{-1}(x)&   \mb{if}~~ x\in\Ima{f} \\
 +\infty && \mb{otherwise}\earr \rp.
\]
For $x \in \Ima{f}$, as $f(y)=f(y')=x$ for each  $y,y' \in f^{-1}(x)$, applying \eq{22} yields $h(y')-h(y) \geq 0$. By symmetry $h(y')=h(y)$, hence the definition of $\ff(x)$ does not depend of which $y \in f^{-1}(x)$ is chosen.  

 For $x'\in\Ima{f} $ we write $x'=f(y')$. Using \eq{22} and the definition of $\ff$ yields 
\[
\ff(x')-\ff(f(y)) = \ff(f(y'))-\ff(f(y)) = h(y')-h(y) \geq 
\<B(f(y'))-B(f(y)),y\> = \<B(x')-B(f(y)),y\>.
\]
that is to say 
\[ \ff(x') - \<B(x'),y\> \geq \ff(f(y)) - \<B(f(y)),y\> ,\qu \all x' \in \Ima{f}.
\]
This also trivially holds for $x' \notin \Ima{f}$. Setting $\ph(x):= \ff(x)-b(x)$ for each $x \in \H$, and  replacing $\ff$ by $b+\ph$ in the inequality above yields
\[
a(y)-\<B(x'),y\>+b(x') + \ph(x') \geq a(y)-\<B(f(y)),y\>+b(f(y))+\ph(f(y)),\qu \all x' \in \H
\]
showing that $f(y) \in \arg\min_{x'} \{{\wt D}(x',y)+\ph(x')\}$. As this holds for each $y \in \Y$, 
$\ph$ satisfies \eq{proabis}.\\

\paragraph{\eq{pro4}}
Consider $\ph$ and $\psi$ satisfying \eq{proabis} and \eq{proc},  respectively. Using the arguments of \eq{proabis} $\Rightarrow$ \eq{proc}, the function $\psi_1 := \breve\gg$ with $\gg$ defined in \eq{psitmp} satisfies \eq{proc}. As $\psi$ and $\psi_1$ both satisfy \eq{proc}, for each $y \in \C'$ we have $\d \psi(y) \cap \d \psi_1(y) \supset B(f(y)) \neq \void$ with $\psi,\psi_1$ convex l.s.c. functions. Hence, by \cor{uni}, since $\C'$ is polygonally connected, there is a constant $K'$ such that $\psi(y) = \psi_1(y)+K'$, $\all y \in \C'$. By \eq{procrhosd}, $\bd \gg(y) \neq \void$ for each $y \in \Y$, hence by \prop{sdcvxenv} we have  $\breve{\gg}(y) = \gg(y)$ for each $y \in \Y$. As $\C' \subset \Y$, it follows that $\psi_1(y) = \breve{\gg}(y) = \gg(y)$ for each $y \in \C'$. This establishes \eq{psi}.

\subsection{Proof of \lem{SdCont}} \label{lsp}
\proof
Without loss of generality we prove the equivalence for the convex envelope $\breve{\ff}$ instead of $\ff$: indeed by \prop{sdcvxenv}, since $\partial \ff(x) \neq \void$ on $\X$ we have $\breve{\ff}(x) = \ff(x)$ and $\partial \breve{\ff}(x) = \partial \ff(x)$ on $\X$.\\
\eq{SdConta}$\Rightarrow$\eq{SdContb}. 
By \cite[Prop 17.41(iii)$\Rightarrow$(i)]{Bauschke:2017ki}, as $\breve{\ff}$ is convex l.s.c. and $\gg$ is a selection of its subdifferential which is continuous at each $x \in \X$, $\breve{\ff}$ is (Fréchet) differentiable at each $x \in \X$. By \prop{sdisg} we get $\partial \breve{\ff}(x) = \{\nabla \breve{\ff}(x)\} = \{\gg(x)\}$ on $\X$. Since $\gg$ is continuous, $x \mapsto \nabla \breve{\ff}(x)$ is continuous on $\X$.\\
\eq{SdContb}$\Rightarrow$\eq{SdConta}. 
Since $\breve{\ff}$ is differentiable on $\X$, by \prop{sdisg} we have $\partial \breve{\ff}(x) = \{\nabla \breve{\ff}(x)\}$ on $\X$. By \eq{ginf} it follows that $\gg(x) = \nabla \breve{\ff}(x)$ on $\X$. Since $\nabla \breve{\ff}$ is continuous on $\X$, so is $\gg$.
\endproof

\subsection{Proof of \cor{uniqueC1}} \lab{pfuniqueC1}

By \thm{jac1}, as $\Y$ is open and convex and $f$ is $C^1(\Y)$ with $Df(y)$ symmetric semi-definite positive for each $y \in \Y$, there is a function $\ph_0$ and a convex l.s.c. function $\psi \in C^2(\Y)$ such that $\nabla\psi(y) = f(y) \in \prox_{\ph_0}(y)$ for each $y \in \Y$. We define $\ph(x) := \ph_0(x) + \chi_{\Ima{f}}(x)$ and let the reader check that $f(y) \in \prox_\ph(y)$ for each $y \in \Y$.  By construction, $\dom{\ph} = \Ima{f}$.

{\bf Uniqueness of the global minimizer.} Consider $\wt{f}$ any function such that $\wt{f}(y) \in \prox_{\ph}(y)$ for each $y$. This implies 
\beq \lab{uC10}
\tfrac{1}{2}\|y-f(y)\|^2+\ph(f(y))
= 
\tfrac{1}{2}\|y-\wt{f}(y)\|^2+\ph(\wt{f}(y))
=
\min_{x \in \H} \{\tfrac{1}{2}\|y-x\|^2+\ph(x)\},\qu \all y \in \Y.
\eeq
By \cor{PCDS} there is a convex l.s.c. function $\wt{\psi}$ such that $\wt{f}(y) \in \d \wt{\psi}(y)$ for each $y \in \Y$. Since $\Y$ is convex it is polygonally connected hence by \thm{ProxCh}\eq{tpro} and \eq{uC10} there are $K,K' \in \RR$ such that
\beq \lab{uC10bis}
\psi(y) -K
= \tfrac{1}{2}\|y\|^2-\tfrac{1}{2}\|y-f(y)\|^2-\ph(f(y)) = \tfrac{1}{2}\|y\|^2-\tfrac{1}{2}\|y-\wt{f}(y)\|^2-\ph(\wt{f}(y)) 
= \wt{\psi}(y) - K',\qu \all y \in \Y.
\eeq
Thus, $\wt{\psi}$ is also $C^2(\Y)$ and $\wt{f}(y) \in \d\wt{\psi}(y) = \{\nabla \psi(y)\} = \{f(y)\}$ for each $y \in \Y$. This shows that $\wt{f}(y)=f(y)$ for each $y$, hence $f(y)$ is the unique global minimizer on $\H$ of $x \mapsto \tfrac{1}{2}\|y-x\|^2+\ph(x)$, i.e., $\prox_\ph(y) = \{f(y)\}$.

{\bf Injectivity of $f$.} The proof follows that of \cite[Lemma 1]{GRIBONVAL:2010:INRIA-00486840:1}. Given $y \neq y'$ define $v := y'-y \neq 0$ and $\ff(t) := \<f(y+tv),v\>$ for $t \in [0,1]$. As $\Y$ is convex this is well defined. As $f \in \C^1(\Y)$ and $Df(y+tv) \succ 0$, the function $\ff$ is $C^1([0,1])$ with $\ff'(t) = \< Df(y+tv)\ v,v\> > 0$ for each $t$. If we had $f(y) = f(y')$ then by Rolle's theorem there would be $t \in [0,1]$ such that $\ff'(t)=0$, contradicting the fact that $\ff'(t)>0$.

{\bf Differentiability of $\ph$.}
If $Df(y)$ is boundedly invertible for each $y \in \Y$, then by the inverse function theorem $\Ima{f}$ is open and $f^{-1}: \Ima{f} \to \Y$ is $C^{1}$. Given $x \in \Ima{f}$, denoting $u := f^{-1}(x)$, \eqref{uC10bis} yields 
\[
\ph(x) = \ph(f(u)) = -(\psi(u)-K)+\tfrac{1}{2}\|u\|^{2}-\tfrac{1}{2}\|u-f(u)\|^{2} = -(\psi(f^{-1}(x))-K)+\tfrac{1}{2}\|f^{-1}(x)\|^{2}-\tfrac{1}{2}\|f^{-1}(x)-x\|^{2}.
\]
Since $\psi$ is $C^{2}$ and $f^{-1}$ is $C^{1}$, it follows that $\ph$ is $C^{1}$.

{\bf Global minimum is the unique critical point.}
The proof is inspired by that of \cite[Theorem 1]{GRIBONVAL:2010:INRIA-00486840:1}. Consider $x$ a critical point of $\ff: x \mapsto \tfrac{1}{2}\|y-x\|^2+\ph(x)$, i.e., since $\ph$ is $C^{1}$, a point where $\nabla \ff(x)=0$. Since $\dom{\ph} = \Ima{f}$ there is some $v \in \Y$ such that  $x = f(v)$. Moreover, as $\ph$ is $C^{1}$ on the open set $\Ima{f}$, the gradient $\nabla \ff(x)$ is well defined and $\nabla \ff(x)=0$. 
On the one hand, denoting $\gg(u):= (\ff \circ f)(u) = \tfrac{1}{2}\|y-f(u)\|^2+\ph(f(u))$ we have $\nabla \gg(u) = Df(u) \nabla \ff(f(u))$ for each $u \in \Y$. On the other hand, for each $u \in \Y$, as $f(u) = \nabla \psi(u)$ we also have
\begin{eqnarray*}
\gg(u) &=& \tfrac{1}{2}\|y\|^2+\tfrac{1}{2}\|f(u)\|^2-\<y,f(u)\>+\ph(f(u))\\
&=&  +\tfrac{1}{2}\|y\|^2+\<u-y,f(u)\>-(\psi(u)-K),\\
\nabla \gg(u) &=& Df(u)\ (u-y) + f(u) -\nabla \psi(u) = Df(u)\ (u-y)
\end{eqnarray*}
For $u=v$ we get $Df(v)\ (v-y) = \nabla \gg(v) = Df(v) \nabla \ff(f(v)) = Df(v)\ \nabla \ff(x) = 0$.
As $Df(v) \succ 0$, this implies $v=y$, hence $x=f(y)$.

\subsection{Proof of \lem{le:PEWNotProx}} \lab{app:PfLePEWNotProx}

As a preliminary let us compute 
the entries of the $n \times n$ matrix associated to $Df(y)$:
\beq \lab{partialPEW} 
\all i,j \in \set{n}\qu 
\tfrac{\partial f_{i}}{\partial y_{j}}(y) =
\begin{cases}
0 & \text{if}\ \|\diag(w^{i})y\|_{2} < \lambda\\
2(w^{i}_{j})^{2} y_{i}y_{j}
h'_{i}\left(\|\diag(w^{i})y\|_2^{2}\right) 
& \text{if}\ \|\diag(w^{i})y\|_{2}>\lambda
\end{cases}
\eeq
NB: if $\|\diag(w^{i})y\|_{2}=\lambda$ then $f$ may not be differentiable at $y$; this case will not be useful below.

The proof exploits \cor{jac2} which shows that if $f$ is a proximity operator then $Df(y)$ is symmetric in each open set where it is well defined.

Let $f$ be a generalized social shrinkage operator as described in \lem{le:PEWNotProx} and consider $\mathcal{G} = \{G_{1},\ldots,G_{p}\}$ the partition of $\set{n}$ into disjoint groups corresponding to the equivalence classes
defined by the equivalence relation between indices: for $i,j \in \set{n}$, $i \sim j$ if and only if $w^i = w^j$. Given $G \in \mathcal{G}$, denote $w^G$ the weight vector shared all $i \in G$. If $f$ is a proximity operator then we show that for each $G \in \mathcal{G}$, we have $\supp(w^{G}) = G$.

For $i \in G$, by \defi{gss} we have $i \in N_i = \supp(w^i) = \supp(w^G)$, establishing that\footnote{The inclusion \eq{PEWtmp1} is true even if $f$ is not a proximity operator.}
\beq\lab{PEWtmp1}
G \subset \supp(w^G).
\eeq

From now on we assume that $f$ is a proximity operator, and consider a group $G \in \mathcal{G}$. To prove that $G = \supp(w^G)$, we will establish that for each $i,j \in \set{n}$ 
\beq \lab{PEWpf1}
\text{if there exists}\ y \in \RR^n\ \text{such that}\ 
\|\diag(w^{j})y\|_{2} \neq \|\diag(w^{i})y\|_{2}
\ \text{then}\ w^i_j = 0\ \text{and}\ w^j_i = 0.
\eeq
To see why it allows to conclude, consider $j \in \supp(w^G)$, and $i \in G$. As $N_i := \supp(w^i) = \supp(w^G)$ we obtain that $j \in N_i$, i.e., $w^i_j \neq 0$. 
By \eq{PEWpf1}, it follows that $ \|\diag(w^{j})y\|_{2} = \|\diag(w^{i})y\|_{2}$ for each $y$. As $w^i,w^j$ have non-negative entries, this means that $w^{i} = w^{j}$. 
As $i \in G$, this implies $j \in G$ by the very definition of $G$ as an equivalence class. This shows $\supp(w^G) \subset G$. Using also \eq{PEWtmp1}, we conclude that $\supp(w^G)=G$.

Let us now prove \eq{PEWpf1}. Consider a given pair $i,j \in \set{n}$. Assume that $\|\diag(w^{j})y\|_{2} \neq \|\diag(w^{i})y\|_{2}$ for at least one vector $y$. Without loss of generality assume that
$a := \|\diag(w^{j})y\|_{2} < \|\diag(w^{i})y\|_{2} =:b$. 
Rescaling $y$ by a factor $c = 2\lambda/(a+b)$ yields the existence of $y$ such that for the considered pair $i,j$
\begin{eqnarray}\label{eq:TmpPEWNotProx}
\|\diag(w^{j})y\|_{2}<\lambda< \|\diag(w^{i})y\|_{2}.
\end{eqnarray}
By continuity, perturbing $y$ if needed we can also assume that for this pair $i,j$ we have $y_{i}y_{j} \neq 0$. 

By~\eq{partialPEW}, as~\eqref{eq:TmpPEWNotProx} holds in a neighborhood of $y$, $f$ is $C^1$ at $y$ and its partial derivatives for the considered pair $i,j$ satisfy
\[ 
\tfrac{\partial f_{i}}{\partial y_{j}}(y) 
= 2(w^{i}_{j})^{2} y_{i}y_{j} h'_{i}\left(\|\diag(w^{i})y\|_2^{2}\right)
\quad \text{and}\quad
\tfrac{\partial f_{j}}{\partial y_{i}}(y) = 0.
\]
Since $f$ is a proximity operator, by \cor{jac2} we have $\tfrac{\partial f_{i}}{\partial y_{j}}(y) = \tfrac{\partial f_{j}}{\partial y_{i}}(y)$. It follows that for the considered pair $i,j$
\[
(w^{i}_{j})^{2} y_{i}y_{j} h'_{i}\left(\|\diag(w^{i})y\|_2^{2}\right) 
= 0.
\]
As $y_iy_j \neq 0$ and $h'_i(t) \neq 0$ for $t \neq 0$, we obtain $w^{i}_{j} = 0$. 

To conclude we now show that $w^j_i = 0$. As $w^i_j=0$, $f_i$ is in fact independent of $y_j$ and $\frac{\partial f_{i}}{\partial y_{j}}$ is {\em identically zero} on $\RR^{n}$. 
By scaling $y$ as needed, we get a vector $y'$ such that $y'_{i}y'_{j} \neq 0$ and 
\[
\lambda < \|\diag(w^{j})y'\|_{2} < \|\diag(w^{i})y'\|_{2}.
\]
Reasoning as above yields $2(w^j_i)^2 y'_j y'_i h'_j\lp(\|\diag(w^j)y'\|_2^2\rp) = \frac{\partial f_{j}}{\partial y_{i}}(y') = \frac{\partial f_{i}}{\partial y_{j}}(y') =0$, hence $w^{j}_{i} = 0$.  
We thus obtain that $w^{i}_{j}=w^{j}_{i}=0$ as claimed, establishing \eq{PEWpf1} and therefore $G = \supp(w^G)$.

\bibliographystyle{plain}
\bibliography{biblio}

\end{document}